\documentclass[11pt]{amsart}
\usepackage{color}
\usepackage{amsmath,amssymb,amsfonts,amsthm,bm}
\usepackage{graphicx}
\usepackage{float}
\usepackage{enumerate}
\usepackage{appendix}
  \usepackage[shortlabels]{enumitem}

\usepackage{color}

\def\R{\mathbb R}

\def\N{\mathbb N}

\def\cC{\mathcal C}
\def\cD{\mathcal D}
\def\cE{\mathcal E}

\def\cP{\mathcal P}

\def\d{\partial}

\renewcommand\div{{\rm div}\,}
\renewcommand\lim{{\rm lim}\,}
\renewcommand\exp{{\rm exp}\,}
\renewcommand\sup{{\rm sup}\,}

\renewcommand\log{{\rm log}\,}

\newcommand{\with}{\quad\!\hbox{with}\!\quad}
\newcommand{\andf}{\quad\!\hbox{and}\!\quad}

\newcommand{\Int}{\displaystyle \int}

\def\dr{\delta\!\rho}
\def\du{\delta\!u}

\def\dw{\delta\!w}

\def\dD{\delta\!D}
\def\dE{\delta\!E}

\def\dP{\delta\!P}

\newtheorem{theorem}{Theorem}[section]

 \theoremstyle{definition}
 
 \theoremstyle{remark}
 \newtheorem{remark}[theorem]{Remark}
 
 \numberwithin{equation}{section}

\newcommand{\wt}{\widetilde}
\newcommand{\eps}{\varepsilon}
\renewcommand{\div}{\mbox{div}\,}
\newcommand{\divv}{\mbox{div}_v\,}

\begin{document}

\title[]{An elementary approach to  the pressureless Euler-Navier-Stokes system}
\author{ Rapha\"el Danchin}
\begin{abstract} 
The  pressureless Euler-Navier-Stokes  system can be obtained formally from the Vlasov-Navier-Stokes system, under the assumption that the 
 distribution function describing the density of particles is monokinetic. Its study has been the subject of several recent papers, which have established the global existence of solutions with high enough regularity, for small initial data. In this work, we demonstrate the global existence of strong solutions in the whole space case, without assuming the initial density to be small and regular: it suffices for it to be bounded and for the total mass to be finite. 
 In passing, we obtain optimal decay estimates for the energy and dissipation functionals. As a corollary, we get a long-time description of the density. All these results are based on an elementary energy method, with no need of sophisticated Fourier analysis tools. 
\end{abstract}
\date\today

\maketitle

In this paper, we address  the question of the global well-posedness and
large-time asymptotics for  the following \emph{pressureless
Euler-Navier-Stokes} system in the whole space $\R^3$:
\begin{equation*}\left\{\begin{aligned}
&\rho_t+\div_{\!x}(\rho  w)=0,\\
&w_t+w\cdot\nabla_{\!x} w=u-w,\\
&u_t+ u\cdot\nabla_{\!x} u-\Delta_{x} u+\nabla_{\!x} P=\rho(w-u),\\
&\div_{\!x} u=0.
\end{aligned}\right. 
\leqno(ENS)
\end{equation*}
The above system is a coupling between the incompressible Navier-Stokes system 
satisfied by the divergence-free vector-field $u=u(t,x)$ 
and the pressure function $P=P(t,x),$ and the pressureless compressible Euler system involving 
the density $\rho=\rho(t,x)\geq0$ and the vector-field $w=w(t,x).$
Here, $t\geq0$ denotes the time variable and $x\in\R^3$ is the space variable. 
The return force corresponding to the right-hand sides of the second
and third equations is the so-called Brinkman force. 
\smallbreak
Studying the above system is motivated by recent works dedicated to the following 
\emph{Vlasov-Navier-Stokes system}: 
$$\left\{\begin{array}{l}
f_t+v\cdot\nabla_xf+\divv\bigl((u-v)f\bigr)=0,\\[1ex]
u_t +\div_{\!x}(u\otimes  u)-\Delta_x u+\nabla_x P=\Int_{\R^3} (v-u)f\,dv,\\[1ex]
\div_{\!x} u=0.\end{array}\right.\leqno(VNS)$$
This system is used  for describing the evolution of  a cloud of non-interacting particles 
that are immersed in an incompressible viscous homogeneous fluid 
(see the detailed derivation of the system in Williams' book \cite{Wi}). 
The kinetic distribution function $f=f(t,x,v)$   represents  the  density of particles with velocity $v\in\R^3,$ that 
are located at $x\in\R^3$ at time $t\geq0,$  while $u=u(t,x)$ stands for the ambient velocity. 
To our knowledge, the mathematical study of (VNS)
has been initiated  with the papers by  Anoshchenko {\it et al} \cite{AB}, and Boudin {\it et al} \cite{BDGM}.  
\smallbreak
Provided  $u$ is sufficiently smooth,  $f$ can be seen as the solution of a linear transport equation, given
by  the formula 
\begin{equation}\label{eq:f}
f(t,x,v)=e^{3t} f_0(Z(0;t,x,v)),\qquad (t,x,v)\in\R_+\times\R^3\times\R^3,\end{equation}
where $Z$ is the flow of the vector-field $F:=(v,u-v)$ in $\R^3\times\R^3$ defined as the solution of the ODE
$$\frac d{ds} Z(s;t,x,v)= F(s, Z(s;t,x,v)),\quad Z(t;t,x,v)=(x,v),\quad (s,t,x,v)\in\R_+\times\R_+\times\R^3\times\R^3.$$
We immediately see from \eqref{eq:f} 
that at some points of the phase space, the distribution $f$ experiences exponential growth with respect to $t.$ 
At the same time, the total mass $\|f(t)\|_{L^1(\R^3\times\R^3)}$ is conserved.  Moreover,  as pointed out by 
the author in \cite{R-VNS} (after prior works by Han-Kwan, Moussa and Moyano \cite{HMM}, and  
Han-Kwan  \cite{DHK}), if we set $\rho_f(t,x):=\int_{\R^3} f(t,x,v)\,dv$ then, in the large time asymptotics,
  the distribution function $f(t,\cdot,\cdot)$ 
is well approximated by the monokinetic distribution $\rho_f(t,\cdot)\otimes\delta_{v=u(t,\cdot)}$ in the sense of some Wasserstein distance. 
However, the discrepancy between $f$ and its monokinetic ansatz decays only algebraically  with respect to time, 
to compare with the exponential concentration that we pointed out before.
Consequently, we can   expect that at `intermediate times',  the distribution function is close to 
a (possibly infinite)  sum of monokinetic states, namely, of the form
$$f(t,x,v)=\int_I \rho^\alpha(t,x)\otimes\delta_{v=w^\alpha(t,x)}\,d\mu(\alpha)$$
for some suitable probability space $(I,\mu).$ The family of pairs $(\rho^\alpha,w^\alpha)$  may be understood as 
representing a multiphasic flow.   
\smallbreak
A work in progress by Baradat, Ertzbischoff and Han-Kwan \cite{BEHK} seeks to rigorously justify this link  with multiphasic flows in the context of
 general Vlasov equations.
The pressureless Euler-Navier-Stokes system under consideration here corresponds to monokinetic 
 distribution functions (meaning that   the set $I$  reduces to a single element) of the Vlasov-Navier-Stokes equations.
In other words,  we postulate that $$f(t,x,v)=\rho(t,x)\,\delta_{v=w(t,x)}.$$ 
The definition of $f$ then guarantees that 
 $$\rho(t,x)=\Int_{\R^3} f(t,x,v)\,dv\andf (\rho w)(t,x)=\Int_{\R^3}  vf(t,x,v)\,dv.$$
 Furthermore, multiplying the first equation  of $(VNS)$ by $v,$ then  integrating  with respect to~$v$ formally gives:
 $$ \partial_t\int_{\R^3} vf\,dv +\div_{\!x}\int_{\R^3} fv\otimes v\,dv +{\int_{\R^3} v\,\divv\bigl((u-v)f\bigr)dv}=0.$$
In the end, integrating  by parts in the last term, we discover that $(\rho,w,u)$ satisfies
{\begin{equation*}\left\{\begin{aligned}
&\rho_t+\div_{\!x}(\rho  w)=0,\\
&(\rho w)_t+\div_{\!x}(\rho w\otimes w)={\rho(u-w)},\\
&u_t+ \div_{\!x}(u\otimes u)-\Delta_{x} u+\nabla_{x} P={\rho(w-u)},\\&\div_{\!x} u=0,
\end{aligned}\right. \leqno(\wt{ENS})
\end{equation*} }
which is nothing other than the conservative form of $(ENS).$
\medbreak
Various hydrodynamic limits of (VNS) have been studied recently. For adequate scalings, the convergence to the Navier-Stokes system (homogeneous or inhomogeneous) or Boussinesq system has been justified in \cite{Ertz,HKM}. 
However, the link between  of (VNS) and (ENS) is of a different nature and, so far, 
exhibiting  a rigorous connection between the two systems has remained an open question.
In this respect, it is worth mentioning Lemarié's  recent work \cite{Lem}, where  particular solutions of (VNS)  are constructed
from some multiphasic extension of~(ENS). 
\medbreak
Here, we  focus on  the Cauchy problem for (ENS) with initial data
 $(\rho_0,w_0,u_0)$ such that $\div u_0=0$ in the case where the fluid domain  is  the whole space $\R^3.$
 This question  has been studied before by Y.-P. Choi et al in \cite{CK,CJK} who 
 established the global existence of strong solutions and algebraic decay rates for $(w,u)$ in the case of  small  $(\rho_0,w_0,u_0)$ in Sobolev spaces $H^s$ with large enough $s,$ 
 with, additionally, $(\rho_0,u_0)$ small in~$L^1.$ This latter condition has been removed 
 by  Huang et al in \cite{HTZ} who also proved some optimality result for the $L^2$ decay of $u.$
 In the same spirit, one can mention the work by Zhai \emph{et al} \cite{Zhai} where data are taken in Besov spaces.
 In \cite{Lem}, Lemari\'e established  similar results for data with critical regularity, improved the aforementioned decay 
 rates and studied the multi-fluid generalization of (ENS) mentioned above.    It should also be noted that System (ENS) is globally well-posed for small  perturbations of $(\rho,w,u)=(\bar\rho,0,0)$ with
 the reference density $\bar\rho$ being a positive constant, see \cite{HTWZ} for more details. 
 In this latter case however, the physical meaning of the connection between (VNS) and (ENS) is unclear since the total mass of the particles
 is infinite.
  \medbreak
 All these results on the Euler-Navier-Stokes system use functional settings far from those for which we know how to solve (VNS). We propose here a framework analogous to that used in \cite{R-VNS}, with a view to rigorously justifying the link between (VNS) and (ENS) in future works.

 To better explain the substance of our results, let us recall the basic conservation laws associated to (ENS). 
 The first one is   the conservation of total mass, namely
\begin{equation}\label{eq:M0}
\|\rho(t)\|_{L^1}=M_0:=\|\rho_0\|_{L^1},\quad\hbox{for all }\ t\geq0.
\end{equation}
The second fundamental one is  the energy balance:
\begin{align}\label{eq:energybalance}
E_0(t)+\int_0^t D_0(\tau)\,d\tau =E_0(0)\with 
&E_0:= \frac12\Bigl( \|\sqrt \rho\, w\|_{L^2}^2+\|u\|_{L^2}^2\Bigr)\\
\andf &D_0:=\|\sqrt\rho\,(w-u)\|_{L^2}^2+ \|\nabla u\|_{L^2}^2.\nonumber
\end{align}
We shall also use the following higher order energy and dissipation functionals:
\begin{equation}\label{eq:E1}
E_1:=\|\sqrt\rho\,(w-u)\|_{L^2}^2+\|\nabla u\|_{L^2}^2\andf 
D_1:=\|\sqrt\rho\,(w-u)\|_{L^2}^2+\|u_t\|_{L^2}^2.
\end{equation}
Our first result states the existence of global strong solutions only assuming that
the data have sufficient integrability at infinity, and that  the initial velocity fields 
$u_0$ and $w_0$ are small enough (in the meaning specified below). 
The initial density, as for it,  can be arbitrarily large and has no smoothness whatsoever : it just 
has  to be in $L^1\cap L^\infty.$ 
\begin{theorem}\label{thm:main1}
 Assume that $\rho_0\in L^1\cap L^\infty$ is nonnegative, that $u_0\in L^{3/2}\cap H^1$ with ${\rm div}\, u_0=0,$  
 that $w_0\in C^{0,1}$ and that  $\sqrt{\rho_0}\, w_0\in L^2.$
There exists a constant $c_0>0$ depending only on $M_0,$ $\|\rho_0\|_{L^\infty}$  and  $\|u_0\|_{L^{3/2}}$ such that if 
\begin{equation}\label{eq:smalldata}
\|u_0\|_{H^1}+\|\sqrt{\rho_0}\,w_0\|_{L^2} +\|w_0\|_{C^{0,1}} \leq c_0
\end{equation}
then System (ENS) supplemented with initial data $(\rho_0,w_0,u_0)$ has a unique global solution  $(\rho,w,u,\nabla P)$  such that:
\begin{itemize}
\item  $\rho\in \cC_b(\R_+;L^1)\cap L^\infty(\R_+\times\R^3),$\smallbreak
\item $w\in \cC_b(\R_+\times\R^3)\cap L^2(\R_+;L^\infty)$ and $\nabla w\in L^\infty(\R_+\times \R^3)\cap L^1(\R_+;L^\infty),$\smallbreak
\item $u\in \cC_b(\R_+;H^1)\cap L^2(\R_+;L^\infty)$ and  $\nabla u\in L^2(\R_+;H^1)\cap L^1(\R_+;L^\infty),$\smallbreak
\item $\nabla P\in L^2(\R_+;L^2).$\end{itemize}
This solution satisfies \eqref{eq:M0}, \eqref{eq:energybalance}, 
and  the following  time decay property:
\begin{equation}\label{eq:energydecay} E_{1}(t)\leq 2e^{E_0(0)}E_1(0)  (1+a_1 t)^{-3/2}\end{equation}
with $a_1$ depending only on $M_0,$ $\|\rho_0\|_{L^\infty},$ $E_0(0),$ $E_1(0)$ and $\|u_0\|_{L^{3/2}}.$
\end{theorem}
\begin{remark} Since the constructed  (Euler) velocity-field $w$ is in  $L^1_{loc}(\R_+;C^{0,1}),$ it has 
a Lipschitz flow $W_t,$ and the density $\rho$ is thus given by the formula
\begin{equation}\label{eq:rhoflow}\rho(t,x)=\rho_0(W_t^{-1}(x))\exp\biggl(\int_0^t (\div w)(s,W_s(W_t^{-1}(x)))\,ds\biggr)\cdotp\end{equation}
This in particular entails that the strict positivity of $\rho$ is preserved during the time evolution. 
Consequently,  if in the previous statement we assume in addition that $\rho_0>0,$ then we are guaranteed that solving $(ENS)$ or its conservative counterpart $(\widetilde{ENS})$ is equivalent.
\end{remark} 
If  $u_0$ has  more integrability,  then one can exhibit the decay rate of $E_0,$ improve the  decay rate of $E_1$ and, 
as a consequence, specify the asymptotic behavior of the density, as stated in the following theorem.
\begin{theorem}\label{thm:main1b}
Let the hypotheses of the previous theorem be in force. 
If, in addition, $u_0$ belongs to $L^1,$ then we also have for some positive real numbers $a_0$ and $a'_1$
depending only on suitable norms of the data:
\begin{equation}\label{eq:energydecaybis} 
E_{0}(t)\leq E_0(0)  (1+a_0 t)^{-3/2}\andf
E_{1}(t)\leq 2e^{E_0(0)}E_1(0)  (1+a'_1 t)^{-5/2}.\end{equation}
Moreover, there exists $\rho_\infty\in L^\infty$ such that $(\rho_\infty-\rho(t))\in  \dot W^{-1,1}\cap \dot H^{-1}$ for all $t\geq0$ and
\begin{equation}\label{eq:convdensity}
\|\rho(t)-\rho_\infty\|_{W^{-1,1}\cap H^{-1}}\leq Ct^{-1/4}\quad\hbox{for all }\ t\geq1.
\end{equation}
\end{theorem}
\begin{remark} As in a number of recent works dedicated to decay 
estimates, the relevant assumption on $u_0$ is a control of low frequencies that may be stated
in terms of homogeneous Besov spaces of the type  $\dot B^{-\sigma}_{2,\infty}$ that correspond to the critical embedding
\begin{equation}\label{eq:embed0}L^p(\R^3)\hookrightarrow \dot B^{-\sigma}_{2,\infty}(\R^3)\with \sigma=\frac3p-\frac 32\andf p\in[1,2).\end{equation}
Such Besov spaces  may be defined from the heat
semi-group $(e^{t\Delta})_{t>0}$ as follows:
$$ \dot B^{-\sigma}_{2,\infty}(\R^3):=\Bigl\{z\in{\mathcal S}'(\R^3)\;:\;
 \|z\|_{\dot B^{-\sigma}_{2,\infty}}:=\underset{t>0}\sup \bigl(t^{\sigma/2}\|e^{t\Delta}z\|_{L^2}\bigr)<\infty\Bigr\}\cdotp$$
In this respect, Theorem \ref{thm:main1} still holds is we assume  $u_0$ to be in $\dot B^{-1/2}_{2,\infty}$
instead of $L^{3/2}$  
while Theorem \ref{thm:main1} with suitably modified decay rates
may be extended to $u_0\in B^{-\sigma}_{2,\infty}$ for any $1<\sigma\leq 3/2.$ In light of \eqref{eq:embed0}, 
the endpoints  $\sigma=3/2$ and $\sigma=1$ 
 correspond to $L^1$ and $L^{6/5},$ respectively. \end{remark}

Let us finally state a global-in-time result pertaining to  the case where $u_0$ has critical 
Besov regularity\footnote{The reader may refer to e.g. \cite[Chap. 2]{BCD} for the definition of general homogeneous
Besov spaces. For the time being, let us just point out that 
$\dot B^{1/2}_{2,1}$ is embedded in the homogeneous Sobolev space $\dot H^{1/2},$ and that these two spaces have the same scaling invariance.}:
\begin{theorem}\label{thm:main2} Let $\rho_0$ and $w_0$ be as above and assume that 
$u_0\in L^{3/2}\cap \dot B^{1/2}_{2,1}.$ There exists a constant $c_0>0$ such that if 
\begin{equation}\label{eq:smalldata2}
E_0(0)+\|u_0\|_{\dot B^{1/2}_{2,1}} +\|w_0\|_{C^{0,1}}\leq c_0,
\end{equation}
then $(ENS)$ has a unique global solution $(\rho,w,u,\nabla P)$ with 
$(\rho,w)$ satisfying the same properties as in Theorem \ref{thm:main1}, 
$$
u\in\cC_b(\R_+; B^{1/2}_{2,1}) \cap L^2(\R_+;\dot B^{3/2}_{2,1})\andf  \nabla u\in L^1(\R_+;L^\infty).$$
Furthermore, \eqref{eq:M0} and \eqref{eq:energybalance} hold, as well as  \eqref{eq:energydecay} for $t\geq1$ 
(and also \eqref{eq:energydecaybis} and \eqref{eq:convdensity} if $u_0\in\dot B^{-3/2}_{2,\infty}$).\end{theorem}
\begin{remark} 
It is known since the pioneering work by Fujita and Kato in \cite{FK}  that  the 
classical incompressible Navier-Stokes equations are well-posed for initial velocity $u_0$
in $\dot H^{1/2}.$   This result has been extended recently to (VNS) by the author in  \cite{R-VNS}.
Whether a Fujita-Kato type result holds true for (ENS) is unclear, though. 
Indeed, if one assumes that   $u_0$ belongs to   $\dot H^{1/2}$ rather than to $\dot B^{1/2}_{2,1},$ then 
one  cannot expect   $\nabla u$ to be in $L^1_{loc}(\R_+;L^\infty)$ as
this property generically fails for the heat or Stokes equations supplemented with data in $\dot H^{1/2}.$ 
However, due to \eqref{eq:rhoflow}, controlling the growth of e.g. $\|\rho(t)\|_{L^\infty}$ 
requires a bound for $\|\div w(t)\|_{L^\infty},$ and thus for $\|\nabla w(t)\|_{L^\infty}$ and $\|\nabla u(t)\|_{L^\infty},$
owing to the coupling in (ENS). 
In the case of (VNS),   a log-Lipschitz control on $u$ (that is true if $u_0$ is in $\dot H^{1/2}$)
suffices to justify the formula \eqref{eq:f}, then to bound the density function associated to $f.$
 \end{remark}
 Let us outline the mean ideas leading to the global existence of a solution 
 satisfying the desired properties.
    The key is to control 
 the following quantities uniformly with respect to~$T$:
 \begin{equation}\label{eq:Lip0}\|\rho\|_{L^\infty(0,T\times\R^3)},\quad \int_0^T\|\nabla u\|_{L^\infty}\,dt\andf
 \int_0^T\|\nabla w\|_{L^\infty}\,dt.\end{equation}
From Relation \eqref{eq:rhoflow}, we gather that the Lipschitz control of $w$ ensures the boundedness of $\rho.$
 Furthermore, since $w$ satisfies a dissipative Burgers equation with source term $u,$   the Lipschitz control of $u$
 gives that of $w.$ 
 Hence, the only difficulty is to bound the middle term of \eqref{eq:Lip0}. 
 In fact, since our final goal is to establish the \emph{global} existence for the \emph{nonlinear} system (ENS), we will need this term to be small, namely
 \begin{equation}\label{eq:Lip}
 \int_0^T\|\nabla u\|_{L^\infty}\,dt \leq \varepsilon\ll1,\qquad T>0.
 \end{equation}
  \medbreak
Achieving  \eqref{eq:Lip}
will be based on 
\begin{itemize}
\item three families of energy estimates involving $u$ and $w$;
\item decay estimates of $E_0$ and $E_1$;
\item interpolation inequalities. 
\end{itemize}
It is noteworthy that the strategy and functional framework used to obtain \eqref{eq:Lip} are very similar to those  we used for (VNS)
 in  \cite{R-VNS} although the system under consideration here is purely hydrodynamic. 
This demonstrates that our approach  is likely the right one for establishing a rigorous link between (VNS) and (ENS).

The guiding principle is to replace $\int_{\R^3} f\,dv$ and  $\int_{\R^3} vf\,dv$ by 
$\rho$ and $\rho w,$ respectively.  
The main  difference is that $w$ has to be treated by hyperbolic techniques while, in \cite{R-VNS}, 
$f$ was seen as the solution of a transport equation with respect to both the space and kinetic variables. 
\medbreak
The rest of the paper unfolds as follows: the next section is dedicated to performing
the key a priori estimates leading to the global existence (and the description of asymptotic behavior) for the first two theorems. 
In Section \ref{s:proofs}, we prove Theorems \ref{thm:main1} and \ref{thm:main1b}: first stability with respect to initial data, and 
uniqueness, then, global existence.
 Section \ref{s:critical} is dedicated to  the critical regularity framework of Theorem \ref{thm:main2}.
 


\section{A priori estimates for the global existence}

We want to establish  a priori estimates for (ENS). 
Having in mind to use these estimates not only for proving the global existence but also for stability and uniqueness issues, 
we  consider in a first time  the following  linear system: 
\begin{equation}\label{eq:LENS} \left\{\begin{aligned}
&\rho_t+\div(\rho\wt w)=0,\\
&w_t+\wt w\cdot\nabla w=z-w+f,\\
&u_t+\wt u\cdot\nabla u-\Delta u+\nabla P=\rho(w-z)+g,\\
&\div u=0,
\end{aligned}\right.
\end{equation}
where the source terms $f,$ $g,$ and vector fields $\wt w,$ $\wt u$ and $z$ are given with $\div \wt u=0.$

\subsubsection*{Step 1: Basic energy identity}

Taking the scalar product of the  $w$  equation with $\rho w$  gives
$$\frac12\frac d{dt} \|\sqrt \rho\, w\|_{L^2}^2+\int \rho w\cdot(\wt w\cdot\nabla  w)\,dx
-\frac12\int \rho_t |w|^2\,dx =\int \rho w\cdot z\,dx -\int \rho|w|^2\,dx+\int \rho f\cdot w\,dx.$$
 Using the first equation of \eqref{eq:LENS}
 and integrating by parts  reveals that the last  two  terms of the left-hand side compensate.
 Hence we have
 $$\frac12\frac d{dt} \|\sqrt \rho\, w\|_{L^2}^2=\int \rho w\cdot z\,dx -\int \rho|w|^2\,dx+\int \rho f\cdot w\,dx.$$
 Similarly, taking the scalar product of the  $u$  equation with $u$ yields
 $$\frac12\frac d{dt}\|u\|_{L^2}^2+\|\nabla u\|_{L^2}^2=-\int \rho z\cdot u\,dx +\int \rho w\cdot u\,dx+\int  g\cdot u\,dx.$$
 Hence, summing these two relations and using the notation defined in \eqref{eq:energybalance}, we get
 \begin{equation}\label{eq:energybasic}
 \frac d{dt}E_0+D_0=\int \rho (u-z)\cdot(u-w)\,dx+\int \rho f\cdot w\,dx+\int g\cdot u\,dx.
\end{equation}

\subsubsection*{Step 2: Energy estimate of the dissipation term}

In this part, we assume that  $z=u$ and $f=g=0$ (we do not need to treat a more general case
for solving (ENS)).   Our aim is to get bounds for the functionals $D_1$ and $E_1$ defined in \eqref{eq:E1}.

 As a start, we take the $L^2$ scalar product of the third equation of \eqref{eq:LENS} with $u_t$ and get:
 $$
 \|u_t\|_{L^2}^2+\frac12\frac d{dt} \|\nabla u\|_{L^2}^2 +\int u_t\cdot(\wt u\cdot \nabla u)\,dx+I=0\with
 I:=\int\rho(u-w)\cdot u_t\,dx.$$
 To handle the term $I,$ we do the following computations that are just based 
 on the use of the $\rho$ and $w$  equations and on integration by parts:
 $$ \begin{aligned}
 I&= \int \rho(u-w)\cdot (u_t-w_t)\,dx +\int \rho(u-w)\cdot w_t\,dx\\
 &=\frac12\frac d{dt}\!\int\rho|u-w|^2dx  -\frac12\int\rho_t|u-w|^2dx +\int\rho|u-w|^2dx
+\int\rho(w-u)\cdot(\wt w\cdot\nabla w)\,dx\\
 &=\frac12\frac d{dt}\!\int\rho|u-w|^2\,dx  +\frac12\int\!\div(\rho \wt w) |w-u|^2dx+\!\int\rho|u-w|^2dx
 +\!\int\rho(w-u)\!\cdot\!(\wt w\!\cdot\!\nabla w)\,dx\\
  &=\frac12\frac d{dt}\!\int\rho|u-w|^2\,dx  +\int\!\rho (w-u)\cdot(\wt w\cdot\nabla(u-w))\,dx+\int\rho|u-w|^2dx\\&\hspace{10cm}
 +\int\rho(w-u)\cdot(\wt w\cdot\nabla w)\,dx\\
  &=\frac12\frac d{dt}\!\int\rho|u-w|^2dx +\int\rho|u-w|^2dx 
   +\int\rho (\wt w\cdot\nabla u)\cdot(w-u)\,dx.\end{aligned}$$
Hence we have 
\begin{multline}\label{eq:H1a}
\frac12\frac d{dt} \Bigl(\|\sqrt\rho\,(w-u)\|_{L^2}^2+\|\nabla u\|_{L^2}^2\Bigr) +
\|u_t\|_{L^2}^2+\|\sqrt\rho\,(w-u)\|_{L^2}^2\\=\int\rho (\wt w\cdot\nabla u)\cdot(u-w)\,dx-\int u_t\cdot(\wt u\cdot \nabla u)\,dx.
\end{multline}	
By combining H\"older and Young inequalities, we see that
$$\begin{aligned}
\int\rho (\wt w\cdot\nabla u)\cdot(u-w)\,dx&\leq\frac12\|\sqrt\rho\,(w-u)\|_{L^2}^2
+\frac12\|\sqrt\rho\,\wt w\|_{L^\infty}^2 \|\nabla u\|_{L^2}^2,\\
-\int u_t\cdot(\wt u\cdot \nabla u)\,dx&\leq \frac12\|u_t\|_{L^2}^2+\frac12\|\wt u\|_{L^\infty}^2 \|\nabla u\|_{L^2}^2.
\end{aligned}$$
Hence, using the notation \eqref{eq:E1}, we conclude that 
\begin{equation}\label{eq:H1b}
\frac d{dt} E_1 + D_1\leq 
\bigl(\|\sqrt\rho\,\wt w\|_{L^\infty}^2+\|\wt u\|_{L^\infty}^2\bigr) \|\nabla u\|_{L^2}^2.
\end{equation}	
 Next, to  estimate $\nabla^2u$ and $\nabla P,$ we  rewrite the equation for $u$ as
 a Stokes system, namely
 \begin{equation}\label{eq:stokes}
 -\Delta u+\nabla P=-u_t-\wt u\cdot\nabla u+\rho(w-u),\qquad \div u=0.
 \end{equation}
  We deduce that 
  $$\begin{aligned}
  \|\nabla^2 u\|_{L^2}^2+\|\nabla P\|_{L^2}^2&= \|u_t+\wt u\cdot\nabla u+\rho(u-w)\|_{L^2}^2\\
  &\leq 3 \|u_t\|_{L^2}^2 +3\|\rho\|_{L^\infty} \|\sqrt\rho\,(u-w)\|_{L^2}^2 +3\|\wt u\|_{L^6}^2\|\nabla u\|_{L^3}^2.
  \end{aligned}
  $$
  Hence, leveraging Young and Gagliardo-Nirenberg inequalities, and the fact that
 $$\|\nabla u\|_{L^3}^2\leq C\|\nabla u\|_{L^2}\|\nabla^2u\|_{L^2},$$ we get
  \begin{equation}\label{eq:H1d}
    \|\nabla^2 u\|_{L^2}^2+2\|\nabla P\|_{L^2}^2\leq 6 \|u_t\|_{L^2}^2 +6\|\rho\|_{L^\infty} \|\sqrt\rho\,(u-w)\|_{L^2}^2 
  +C\|\wt u\|_{L^6}^4\|\nabla u\|_{L^2}^2.\end{equation}
  Consequently, assuming that we have for some $R\geq1,$
  \begin{equation}\label{eq:R}
  \|\rho\|_{L^\infty}\leq R,\end{equation}
  we obtain that
  $$
  \frac{\|\nabla^2 u\|_{L^2}^2}{24R} +   \frac{\|\nabla P\|_{L^2}^2}{12R}\leq \frac{\|u_t\|_{L^2}^2}4+\frac{\|\sqrt\rho\,(u-w)\|_{L^2}^2}{4}+
  \frac CR\|\wt u\|_{L^6}^4\|\nabla u\|_{L^2}^2.$$
 Combining with \eqref{eq:H1b}, we conclude that 
 \begin{multline*}
\frac d{dt} E_1 + \wt D_1\leq  \bigl(\|\sqrt\rho\,\wt w\|_{L^\infty}^2+\|\wt u\|_{L^\infty}^2+CR^{-1}\|\wt u\|_{L^6}^4\bigr) \|\nabla u\|_{L^2}^2
\\\with \wt D_1:=\frac{D_1}2+  \frac{\|\nabla^2 u\|_{L^2}^2}{24R} +   \frac{\|\nabla P\|_{L^2}^2}{12R}\cdotp
 \end{multline*}
 For (ENS) (i.e. $\wt w=w$ and $\wt u=u$),  we thus have, using the embedding $\dot H^1\hookrightarrow L^6,$ 
 \begin{equation}\label{eq:H1e}
 \frac d{dt} E_1 + \wt D_1\leq  \bigl(\|\sqrt\rho\,w\|_{L^\infty}^2+\|u\|_{L^\infty}^2\bigr)\|\nabla u\|_{L^2}^2
  +CR^{-1}\|\nabla u\|_{L^2}^6.
 \end{equation}
   Note that one can  estimate $\sqrt \rho\, w_t$ by means of the following relation:
  $$
  \sqrt\rho\, w_t=\sqrt\rho\,(u-w)-\sqrt\rho\, \wt w\cdot\nabla w,$$
   which leads if $\wt w=w$  to 
     \begin{align}\label{eq:H1f}\nonumber
  \|\sqrt\rho\, w_t\|_{L^2}^2&\leq 2\|\sqrt\rho\,(u-w)\|_{L^2}^2+2\|\nabla w\|_{L^\infty}\|\sqrt \rho \,w\|_{L^2}^2\\
  &\leq 2D_1+4\|\nabla w\|_{L^\infty} E_0.\end{align}

\subsubsection*{Step 3:  A higher order  time weighted energy inequality}
 
 From now on, we focus on $(ENS)$ (that is, we consider \eqref{eq:LENS} with $\wt w=w,$  $z=\wt u=u$ and $f=g=0$).  We assume that
 we are given a smooth and sufficiently decaying solution $(\rho,w,u,\nabla P)$ on $[0,T]\times\R^3.$ 
 
 Our aim is to go one step further in the energy method, in order to eventually get a pointwise control of the $H^2$ norm of $u.$
However, since  $u_0$ is only in $H^1,$ we have to compensate this extra degree of regularity by 
multiplying by the weight $\sqrt t,$ owing to the parabolic nature of the equation for $u$. 

With this heuristic in mind, we differentiate   the second and third equations of (ENS) with respect to time, then  multiply by $\sqrt t,$ to get:
\begin{equation*}\left\{\begin{aligned}
&\d_t(\sqrt t\,w_t)=\sqrt t\,(u_t-w_t)-\sqrt t\, w_t\cdot\nabla w-\sqrt t\, w\cdot\nabla w_t +\frac{1}{2\sqrt t} w_t,\\
&\d_t(\sqrt t\, u_t)+  u\cdot\nabla \sqrt t\, u_t-\Delta \sqrt t\, u_t+\nabla \sqrt t\, P_t +\rho\sqrt t\,(u_t-w_t)\\
&\hspace{7cm}=\sqrt t\, (w-u) \rho_t - \sqrt t\, u_t\cdot\nabla u +\frac{1}{2\sqrt t} u_t.
\end{aligned}\right. 
\end{equation*}
Hence, using $\rho_t=-\div(\rho w),$ taking the scalar product of the first equation with $\rho \sqrt t\,w_t,$ of the second one with $\sqrt t\,u_t,$ and summing up, 
\begin{multline}\label{eq:ut}
\frac12\frac d{dt}\Bigl(\|\sqrt t\, u_t\|_{L^2}^2+\|\sqrt{\rho t}\, w_t\|_{L^2}^2\Bigr) +\|\nabla \sqrt t\, u_t\|_{L^2}^2
+\|\sqrt{\rho t}\, (u_t-w_t)\|_{L^2}^2 \\= \frac12\Bigl(\|u_t\|_{L^2}^2+\|\sqrt{\rho}\, w_t\|_{L^2}^2\Bigr) 
-\int t(u_t\cdot\nabla u)\cdot u_t\,dx +\int t(u-w)\cdot u_t\,\div(\rho w)\,dx \\-\int t(w_t\cdot\nabla w)\cdot\rho w_t\,dx
-\int t(w\cdot\nabla w_t)\cdot \rho w_t\,dx-\frac12 \int t\div(\rho w)\,|w_t|^2\,dx.
\end{multline}
We see from integration by parts that the last two terms  compensate each other. 
Next, it is obvious that 
\begin{align}\label{eq:ut1}
-\int t(w_t\cdot\nabla w)\cdot\rho w_t\,dx&\leq \|\nabla w\|_{L^\infty} \|\sqrt{\rho t}\, w_t\|_{L^2}^2\\
-\int t(u_t\cdot\nabla u)\cdot u_t\,dx&\leq \|\nabla u\|_{L^\infty} \|\sqrt{t}\, u_t\|_{L^2}^2,
\label{eq:ut2}\end{align}
and we have (integrate by parts) that
$$\displaylines{\int t(u-w)\cdot u_t\,\div(\rho w)\,dx=\int \rho t (w\cdot\nabla u_t)\cdot (w-u)\,dx
\hfill\cr\hfill+\int t u_t\cdot(\rho (w-u)\cdot\nabla(w-u))\,dx+\int t u_t\cdot(\rho u\cdot\nabla(w-u))\,dx.}$$
On the one hand, by H\"older and Young inequality, we have 
 \begin{align}\nonumber\label{eq:ut3}
\int \rho t (w\cdot\nabla u_t)\cdot (w-u)\,dx&\leq \|\sqrt t\,\nabla u_t\|_{L^2} \|\sqrt{\rho t}\, (w-u)\|_{L^2}\|\sqrt \rho \, w\|_{L^\infty}\\
&\leq\frac12\|\sqrt t\,\nabla u_t\|_{L^2}^2 +\|\sqrt \rho\, w\|_{L^\infty}^2
\|\sqrt{\rho t}(w-u)\|_{L^2}^2.
\end{align}
On the other hand, using also Sobolev inequality, 
\begin{align}\label{eq:ut4}
\int t u_t\cdot(\rho (w-u)\cdot\nabla(w-u))\,dx&\leq \|\rho\|_{L^{\infty}}^{1/2} \|\nabla(w-u)\|_{L^\infty}\|\sqrt t\, u_t\|_{L^2} 
\|\sqrt{\rho t}\,(w-u)\|_{L^2},\end{align}
and finally, 
 \begin{align}\nonumber\label{eq:ut5}
\int t u_t\cdot(\rho u\cdot\nabla(w-u))\,dx&\leq \|\rho\|_{L^3} \|\sqrt t\,u\|_{L^6}  \|\sqrt t\,u_t\|_{L^2} \|\nabla(w-u)\|_{L^\infty}\\
&\leq C \|\rho\|_{L^3} \|\nabla(w-u)\|_{L^\infty}\|\sqrt t\nabla u\|_{L^2}\|\sqrt u_t\|_{L^2}.
\end{align}
Hence,  denoting 
$$E_2:=\|u_t\|_{L^2}^2+\|\sqrt{\rho}\, w_t\|_{L^2}^2\andf 
D_2:=\|\nabla u_t\|_{L^2}^2 +2\|\sqrt{\rho}\, (u_t-w_t)\|_{L^2}^2,$$
and inserting \eqref{eq:ut1} to \eqref{eq:ut5} in \eqref{eq:ut}, we end up with 
\begin{multline}\label{eq:ut6}
\frac d{dt}(tE_2) + tD_2\leq \|u_t\|_{L^2}^2+\|\sqrt{\rho}\, w_t\|_{L^2}^2+
2\|\nabla w\|_{L^\infty} \|\sqrt{\rho t}\, w_t\|_{L^2}^2+2\|\nabla u\|_{L^\infty} \|\sqrt{t}\, u_t\|_{L^2}^2
\\+2\|\sqrt \rho\, w\|_{L^\infty}^2\|\sqrt{\rho t}(w-u)\|_{L^2}^2+
 \|\rho\|_{L^{\infty}}^{1/2} \|\nabla(w-u)\|_{L^\infty}\|\sqrt t\, u_t\|_{L^2} 
\|\sqrt{\rho t}\,(w-u)\|_{L^2}\\
+ C \|\rho\|_{L^3} \|\nabla(w-u)\|_{L^\infty}\|\sqrt t\nabla u\|_{L^2}\|\sqrt t\,u_t\|_{L^2}.
\end{multline}
Since $R\geq1$ and $\|\rho\|_{L^3}\leq \|\rho\|_{L^1}^{1/3}\|\rho\|_{L^\infty}^{2/3},$ using \eqref{eq:M0}, \eqref{eq:R},  and  the definition of $D_1,E_1,$ this inequality may be rewritten
\begin{multline}\label{eq:ut7}
\frac d{dt}(tE_2) + tD_2\leq  D_1+\|\sqrt{\rho}\, w_t\|_{L^2}^2+
2t\|(\nabla w,\nabla u)\|_{L^\infty} E_2+2tR \|w\|_{L^\infty}^2E_1 \\+Ct(R+M_0)\|\nabla(w-u)\|_{L^\infty}\sqrt{E_1 E_2}.
\end{multline}
Putting together \eqref{eq:energybasic} with $z=\wt u=u$ and $\wt w=w,$   \eqref{eq:H1b}  and 
\eqref{eq:ut7}, we conclude that
\begin{multline}\label{eq:cE}
\frac d{dt}\cE+\cD \leq C\Bigl(R\|(w,u)\|_{L^\infty}^2+(R+M_0)\|(\nabla w,\nabla u)\|_{L^\infty}\Bigr)\cE+CR^{-1}D_0\cE^2\\
\with \cE:=2E_0+6E_1 + t (E_1+E_2)\andf \cD:= D_0+D_1+\|\sqrt\rho\, w_t\|_{L^2}^2 +t(D_1+D_2).
\end{multline}
It is now clear that if $\cE(0)$ is small enough, then 
Gronwall lemma enables us to bound the solution on $[0,T]$  in terms of the initial data, of $R$ 
and of the following quantity:
\begin{equation}\label{eq:control}
A_T:=\|(\nabla w,\nabla u)\|_{L^1(0,T;L^\infty)}+\|(w,u)\|_{L^2(0,T;L^\infty)}^2.\end{equation}

\subsubsection*{Step 4: Bounding $A_T$}

Provided \eqref{eq:R} holds true, bounding $u$ in $L^2(0,T;L^\infty)$ 
follows from the previous steps and the following Gagliardo-Nirenberg inequality: 
\begin{equation}\label{eq:GN}
\|u\|_{L^\infty}^2
 \lesssim \|\nabla u\|_{L^2}\|\nabla^2u\|_{L^2}.\end{equation}
We readily get
\begin{equation}\label{eq:vL2L3}
\|u\|_{L^2(0,T;L^\infty)}^2
\leq C\int_0^T\sqrt{RD_0\wt D_1}\,dt.
\end{equation}
The last term may be bounded by means of  \eqref{eq:H1e} provided we have \eqref{eq:R} and $E_1(0)$ is small enough.
Now, the key to controlling the other terms in \eqref{eq:control} is to ensure 
 that $\nabla w$ is small enough in $L^1(0,T;L^\infty).$ This will follow from the fact that 
 the same property holds for $\nabla u.$ 
 Indeed,  the equation for $w$ may be written as
$$\d_t(e^t w)+w\cdot\nabla(e^t w) = e^tu,$$
whence, for all $t\in[0,T],$ 
$$
\|e^t w(t)\|_{L^\infty}\ \leq \|w_0\|_{L^\infty}  +\int_0^t e^\tau\|u\|_{L^\infty}\,d\tau
+\int_0^t \|\nabla w\|_{L^\infty}\|e^\tau w\|_{L^\infty}\,d\tau.
$$
By Gronwall lemma, this gives 
\begin{equation}\label{eq:uL}
\|w(t)\|_{L^\infty}\leq e^{-t}\|w_0\|_{L^\infty}
+\int_0^t e^{\tau-t} \|u(\tau)\|_{L^\infty}\,d\tau.
\end{equation}
Therefore we have, taking advantage of some convolution inequality, 
\begin{equation}\label{eq:uL2}
\|w\|_{L^2(0,T;L^\infty)} + \|w\|_{L^\infty(0,T;L^\infty)} \leq \bigl(\|w_0\|_{L^\infty}+\|u\|_{L^2(0,T;L^\infty)}\bigr)\cdotp 
\end{equation}
Similarly, after differentiating once the equation of $w,$ we  get for all  $t\in[0,T],$
\begin{equation}\label{eq:dwLinfty}e^t\bigl\|\nabla w(t)\bigr\|_{L^\infty}\leq  \bigl\|\nabla w_0\bigr\|_{L^\infty}
+\int_0^te^{\tau}  \bigl\|\nabla u\bigr\|_{L^\infty}
+\int_0^t  \bigl\|\nabla w\bigr\|_{L^\infty} \: e^\tau   \bigl\|\nabla w\bigr\|_{L^\infty}\,d\tau.\end{equation}
Hence, combining Gronwall lemma and convolution inequalities yields
\begin{multline}
\label{eq:duinfty}
\max\Bigl(\|\nabla w\|_{L^\infty(0,T; L^\infty)}, \|\nabla w\|_{L^1(0,T;L^\infty)}\Bigr)\\\leq \bigl(\|\nabla w_0\|_{L^\infty}+ \|\nabla u\|_{L^1(0,T; L^\infty)}\bigr)e^{\|\nabla w\|_{L^1(0,T;L^\infty)}}.
\end{multline}
Assume that
\begin{equation}\label{eq:eps}
2\int_0^T\!\!\|\nabla w\|_{L^\infty}\,dt\leq \log2,\quad 
\int_0^T\!\!\|\nabla u\|_{L^\infty}\,dt\leq \eps\andf \|\nabla w_0\|_{L^\infty}\leq\eps
\quad\!\!\!\!\hbox{for some}\ \eps>0.
\end{equation}
Then, Inequality \eqref{eq:duinfty} guarantees that 
$$\|\nabla w\|_{L^1(0,T;L^\infty)}\leq 4\eps,$$
and provided $\eps$ is chosen so that $8\eps<\log 2,$ we have the first part  of \eqref{eq:eps} with strict inequality. 
Hence, assuming that $w_0$ satisfies the third inequality of \eqref{eq:eps}, it is only a matter of  proving  the middle one.
Finally, from the density equation,  we have  for all $t\geq0$
\begin{equation}\label{eq:rhoLinfty}
\|\rho(t)\|_{L^\infty}\leq \|\rho_0\|_{L^\infty} +\int_0^t\|\div w\|_{L^\infty}\|\rho\|_{L^\infty}\,d\tau.
\end{equation}
Therefore, the above Lipschitz control of $w$ 
ensures \eqref{eq:R}  with $R=\max(1,2\|\rho_0\|_{L^\infty}).$
\smallbreak
In a nutshell, in order to bound all the terms in  \eqref{eq:control}, it is enough to ensure  that 
$u$ satisfies \eqref{eq:Lip}. 
However, this cannot be achieved just by combining  the previous bounds with suitable functional embedding: 
we need to glean more time decay for $u.$ This is the purpose of the next step.

\subsubsection*{Step 5: Time decay}
As in \cite{R-VNS},  we first show that one can propagate some negative Besov regularity, leveraging
 maximal regularity estimates in Besov spaces  for the heat equation:
 from e.g. \cite{BCD,Chemin}, we know that for all $t\in[0,T],$
\begin{equation}\label{eq:v12}\|u\|_{L^\infty(0,t;\dot B^{-1/2}_{2,\infty})} \leq \|u_0\|_{\dot B^{-1/2}_{2,\infty}} +  \|u\cdot\nabla u\|_{L^2(0,t;\dot B^{-3/2}_{2,\infty})}
+\|\rho(w-u)\|_{L^2(0,t;\dot B^{-3/2}_{2,\infty})}.\end{equation} 
Now, owing to the embedding \eqref{eq:embed0}
 with $p=1$ and H\"older inequality, we have
$$\|u\cdot\nabla u\|_{L^2(0,t;\dot B^{-3/2}_{2,\infty})}\lesssim \|u\cdot\nabla u\|_{L^2(0,t;L^1)}
\leq \|u\|_{L^\infty(0,t;L^2)}\|\nabla u\|_{L^2(0,t;L^2)}\leq E_0(0).$$
In light of  \eqref{eq:M0} and \eqref{eq:energybasic}, we also have 
$$\begin{aligned}\|\rho(w-u)\|_{L^2(0,t;\dot B^{-3/2}_{2,\infty})}&\lesssim \|\rho(w-u)\|_{L^2(0,t;L^1)}\\
&\leq \|\sqrt\rho\|_{L^\infty(0,t;L^2)} \|\sqrt \rho(w-u)\|_{L^2(0,t;L^2)}\leq \sqrt{M_0 E_0(0)}.\end{aligned}$$
Hence, reverting to \eqref{eq:v12} gives for some absolute constant $C$ and all $t\in[0,T],$
\begin{equation}\label{eq:v323}
\|u(t)\|_{\dot B^{-1/2}_{2,\infty}} \leq C_0:=  C\bigl(\|u_0\|_{\dot B^{-1/2}_{2,\infty}} + \sqrt{E_0(0)}\bigl(\sqrt{M_0}+\sqrt{E_0(0)}
\bigr)\bigr)\cdotp
\end{equation}
This control of the negative Besov norm readily supplies some time decay for $E_1.$ 
Indeed,  assuming that $E_1$ is small enough, Inequality  \eqref{eq:H1e} implies that 
$$\frac d{dt}E_1+\wt D_1\leq CU E_1\with U:=\|\sqrt\rho\, w\|_{L^\infty}^2+ \|u\|_{L^\infty}^2+D_0.$$
Then, we have
\begin{equation}\label{eq:wtE1}
\frac d{dt}\check E_1+\check D_1\leq 0\with \check E_1(t):= e^{-C\int_0^t U\,d\tau}E_1(t)\andf  \check D_1(t):= e^{-C\int_0^t U\,d\tau}\wt D_1(t).\end{equation}
In the following computations, we make the following a priori assumption: 
\begin{equation}\label{eq:U}
C\int_0^T\bigl(\|\sqrt \rho\,w\|_{L^\infty}^2 + \|u\|_{L^\infty}^2\bigr)\,d\tau\leq \log2.
\end{equation}
Recall the interpolation inequality
\begin{equation}\label{eq:interpo}
\|z\|_{L^2}\lesssim \|\nabla z\|_{L^2}^{3/5}\|z\|_{\dot B^{-3/2}_{2,\infty}}^{2/5}.\end{equation}
Owing to \eqref{eq:v323} and to \eqref{eq:interpo} with $z=\nabla u,$ we thus have after performing  a harmless change of $C_0$: 
$$\|\nabla ^2 u\|_{L^2}^2\geq C_0^{-4/3}\bigl(\|\nabla u\|_{L^2}^2\bigr)^{5/3}.$$ 
Reverting to \eqref{eq:wtE1} and using \eqref{eq:energybalance},  we get
$$\frac d{dt}\check E_1+\frac{e^{-E_0(0)}}{2}\Bigl(C_0^{-4/3}\bigl(\|\nabla u\|_{L^2}^2\bigr)^{5/3}+\|\sqrt\rho\,(w-u)\|_{L^2}^2\Bigr)\leq0.$$
Since $\wt E_{1}$ is nonincreasing and \eqref{eq:U} holds, we have
$$\|\sqrt\rho\,(w-u)\|_{L^2}^2\leq2e^{E_0(0)}\check E_1 \leq 2 e^{E_0(0)}E_1(0).$$ 
Hence if, say, 
\begin{equation}\label{eq:E10} E_1(0)\leq e^{-E_0(0)}/2,\end{equation}
then 
$$1\geq \|\sqrt\rho\,(w-u)\|_{L^2}^2\geq \bigl(\|\sqrt\rho\,(w-u)\|_{L^2}^2\bigr)^{5/3},$$ and thus
$$\frac d{dt}\check E_1+ c_0(\check E_1)^{5/3}\leq 0\with c_0:=2^{-5/3}e^{-E_0(0)}\min(1,C_0^{-4/3}).$$ 
One can conclude that whenever \eqref{eq:U} and \eqref{eq:E10} hold,  we have
\begin{equation}\label{eq:decay}E_1(t)\leq 
2e^{E_0(0)}\bigl(1+a_0 t\bigr)^{-3/2}  E_1(0)\with a_0:=\frac13 c_0  E_1^{2/3}(0).\end{equation}
Note that \eqref{eq:wtE1} implies that for all $\beta>0,$
$$\frac d{dt}(1+a_0t)^\beta  \check E_1+ (1+a_0t)^\beta \check D_1\leq \beta a_0 (1+a_0t)^{\beta-1} \check E_1.$$
Hence, taking $\beta\in(0,3/2)$ and integrating on $[0,T]$ we discover that 
\begin{equation}\label{eq:decay+}
\int_0^T   (1+a_0t)^\beta \check D_1(t)\,dt \leq \frac{\beta}{3/2-\beta} \: E_1(0).\end{equation}

\subsubsection*{Step 6: The Lipschitz bound}
We now have enough information on the solution to prove  \eqref{eq:Lip} in the case of sufficiently small $E_1(0).$ 
Remembering  the critical embedding inequality
\begin{equation}\label{eq:embed} \|z\|_{L^\infty}\leq C\|\nabla z\|_{L^{3,1}},\end{equation}
where  $L^{3,1}$ stands for the  Lorentz space defined by real interpolation as follows:
\begin{equation}\label{def:L31}
L^{3,1}:=[L^2,L^6]_{1/2,1},\end{equation}
it suffices to establish that 
\begin{equation}\label{eq:Lip31}
C\int_0^T\|\nabla^2 u\|_{L^{3,1}}\leq \eps.
\end{equation}
To do this, we use the elliptic maximal regularity for  the Stokes system \eqref{eq:stokes} in $L^{3,1}$
(a straightforward consequence of the classical result in Lebesgue spaces and of the real interpolation theory)
 getting:
\begin{equation}\label{eq:stokes31}
\|\nabla^2u,\nabla P\|_{L^{3,1}} \lesssim \|u\cdot\nabla u\|_{L^{3,1}} + \|u_t\|_{L^{3,1}} + \|\rho(w-u)\|_{L^{3,1}}.\end{equation}
In order to estimate the first term, we use H\"older inequality and \eqref{eq:embed} as follows:
$$
 \|u\cdot\nabla u\|_{L^{3,1}} \leq \|u\|_{L^\infty}\|\nabla u\|_{L^{3,1}}\lesssim \|\nabla u\|_{L^{3,1}}^2.
 $$
 Now, we have the following interpolation inequality~: 
 \begin{equation}\label{eq:L31}\|z\|_{L^{3,1}}^2\lesssim \|z\|_{L^2}\|\nabla z\|_{L^2},\end{equation}
 and thus, denoting  $\wt E_2:= E_2+\|\nabla^2 u\|_{L^2}^2$ and using the definition of $\wt D_1,$
  we get\footnote{In all this paragraph, we agree that integrals on $[1,T]$ are null if $T\leq1.$}:
 $$\begin{aligned}
 \int_0^T  \|u\cdot\nabla u\|_{L^{3,1}}\,dt&\lesssim \int_0^1 \|\nabla u\|_{L^2}\|\sqrt t\, \nabla^2u\|_{L^2} t^{-1/2}dt
 + \int_1^T \|\nabla u\|_{L^2}\|\sqrt t\, \nabla^2u\|_{L^2} t^{-1/2}dt\\
 &\lesssim  \|\nabla u\|_{L^\infty(L^2)}\|\sqrt t\, \nabla^2u\|_{L^\infty(L^2)} 
 + \|\nabla u\|_{L^2(L^2)}\|\sqrt t\, \nabla^2u\|_{L^2(L^2)}\\
 &\lesssim \sqrt{\underset{t\in(0,1)}\sup E_1(t) \: \underset{t\in(0,1)}\sup t\wt E_2(t)} +\biggl(\int_1^T D_0\,d\tau
 \int_1^T tR\wt D_1\,d\tau\biggr)^{1/2}.\end{aligned}$$
According to \eqref{eq:H1d} and to the definition of $E_1$ and $E_2,$ we have
$$\|\nabla^2 u\|_{L^2}^2\leq 6E_2+6\|\rho\|_{L^\infty} E_1+CE_1^3.$$
Hence, provided all the norms in  \eqref{eq:control} are, say, smaller than $1,$ 
Inequality \eqref{eq:cE} guarantees that there exists some constant $K_0$ depending only on the norms
of the data and tending to $0$ when $E_1(0)$ goes to $0,$ such that
\begin{equation}\label{eq:vnablav} 
 \int_0^T  \|u\cdot\nabla u\|_{L^{3,1}}\,dt\leq K_0.\end{equation}
 Next, by definition of the Lorentz space $L^{3,1},$ we have
$$ \|\rho(w-u)\|_{L^{3,1}}\lesssim   \|\rho(w-u)\|_{L^{2}}^{2/3}\|\rho(w-u)\|_{L^{\infty}}^{1/3}.$$
Therefore,  by H\"older inequality,  we have for all $\beta>0$ and $\eps>0,$ 
\begin{align}
\int_0^T  \|\rho(w-u)\|_{L^{3,1}}\,dt&\lesssim \|\rho\|_{L^\infty}^{2/3}\biggl(\int_0^1 \|\sqrt\rho\,(w-u)\|_{L^{2}}^{2/3}\|(w,u)\|_{L^{\infty}}^{1/3}\,dt
\nonumber\\&\hspace{3cm}+
\int_1^T \|t^\beta\sqrt{\rho}\,(w-u)\|_{L^{2}}^{2/3}\|t^\eps(w,u)\|_{L^{\infty}}^{1/3}\,t^{-\frac{\eps+2\beta}3}dt\biggr)\nonumber\\\nonumber
&\lesssim \|\rho\|_{L^\infty}^{2/3} \biggl(\|\sqrt\rho(w-u)\|_{L^2(0,1;L^2)}^{2/3} \|(w,u)\|_{L^2(0,1;L^\infty)}^{1/3}\|1\|_{L^2(0,1)}^{1/2}
\\&\hspace{0.7cm}+\|t^\beta\sqrt{\rho}\,(w-u)\|_{L^2(1,T;L^2)}^{2/3} \|t^\eps(w,u)\|_{L^2(1,T;L^\infty)}^{1/3}\|t^{-\frac{\eps+2\beta}3}\|_{L^{2}(1,T)}\biggr)\cdotp\label{eq:rhowu}
\end{align}
Note that $\|t^{-\frac{\eps+2\beta}3}\|_{L^{2}(1,T)}<\infty$ whenever $2\beta+\eps>3/2$ and 
that $\|t^\beta\sqrt{\rho}\,(w-u)\|_{L^2(1,T;L^2)}$ may be bounded by means of \eqref{eq:decay+} 
for any $\beta<3/4.$ Hence, in order to conclude, it suffices to find some $\eps>0$ such that
both $t^\eps w$ and $t^\eps u$ are in $L^2(1,T;L^\infty).$ 
Now, by Gagliardo-Nirenberg inequality \eqref{eq:GN}, we have 
for all $\eps\in(0,3/8),$ 
$$\begin{aligned}\int_1^T (t^\eps\|u\|_{L^\infty})^2\,dt 
&\lesssim \int_1^T \|t^{3/4}\nabla u\|_{L^2} \|t^{1/2}\nabla^2u\|_{L^2} t^{2\eps-5/4}\,dt\\
&\lesssim \biggl(\underset{t\in(0,T)}\sup t^{3/4} \|\nabla u(t)\|_{L^2}\biggr)  \biggl(\int_0^T \|t^{1/2}\nabla^2u\|_{L^2}^2\,dt\biggr)^{1/2}\\
&\lesssim \biggl(\underset{t\in(0,T)}\sup t^{3/2} E_1(t) \biggr) \biggl(\int_0^T R\wt D_1\,dt\biggr)^{1/2}\cdotp
  \end{aligned}$$
The right-hand side may be bounded in terms of the initial data, owing to our decay estimate \eqref{eq:decay} and to \eqref{eq:H1e}.
Next, to bound $t^\eps w$ in  $L^2(\R_+;L^\infty),$ it suffices to observe that
$$\d_t(t^\eps w)+w\cdot\nabla(t^\eps w) + t^\eps w=\eps w t^{\eps-1}+t^\eps u,$$
which implies that for all $t\geq1,$
$$\displaylines{\quad \|t^\eps w(t)\|_{L^\infty}\leq   e^{\int_1^t \|\nabla w\|_{L^\infty}\,d\tau} e^{1-t}\biggl( \|w(1)\|_{L^\infty} 
\hfill\cr\hfill+ \int_1^t e^{-\int_1^\tau\|\nabla w\|_{L^\infty}\,d\tau'} e^{\tau-1}\Bigl(\|t^\eps u\|_{L^\infty} 
+\eps \tau^{\eps-1}\|w\|_{L^\infty}\Bigr)\,d\tau\biggr)\cdotp\quad}$$
Hence, assuming  \eqref{eq:eps}, we end up with
$$\|t^\eps w\|_{L^2(0,t;L^\infty)}\leq 2\Bigl(\|w(1)\|_{L^\infty} +  \|t^\eps u\|_{L^2(1,t;L^\infty)}+ \varepsilon\|w\|_{L^2(1,t;L^\infty)}\Bigr)\cdotp$$
The right-hand side may be bounded thanks to the above bound (for $u$) 
and to   \eqref{eq:uL2}, hence   in terms of the data. 
Reverting to \eqref{eq:rhowu}, one can conclude that 
 \begin{equation}\label{eq:rhouvL31} \int_0^T  \|\rho(w-u)\|_{L^{3,1}}\,dt\leq K_0.\end{equation}
 
Let us finally bound the term of \eqref{eq:stokes31} pertaining to $u_t.$ Using again \eqref{eq:L31} and H\"older inequality, one can write
for all $\beta\in(1,3/2),$ 
  $$\begin{aligned}
 \int_0^T \! \|u_t\|_{L^{3,1}}\,dt&\lesssim \int_0^1\! \|\sqrt t\,u_t\|_{L^2}^{1/2}\|\sqrt t\,\nabla u_t\|_{L^2}^{1/2} \,t^{-1/2}dt
 \!+\!\int_1^T \!\|t^{\beta/2} u_t\|_{L^2}^{1/2} \|\sqrt t\,\nabla u_t\|_{L^2}^{1/2}\,t^{-(\beta+1)/4} dt\\
 &\lesssim \biggl(\Bigl(\underset{t\in[0,1]}\sup tE_2(t)\Bigr) \int_0^1  tD_2(t)\,dt\biggr)^{1/4} +C_\beta \biggl(\int_1^T t^\beta D_1(t)\,dt \int_1^T  tD_2(t)\,dt\biggr)^{1/4},
\end{aligned}$$
where we used the fact that $t^{-1/2}\in L^{4/3}(0,1)$ and that 
 $t^{-(\beta+1)/4}$ is in $L^2(1,T)$ for $\beta>1.$
 \medbreak
 Remembering \eqref{eq:cE} and \eqref{eq:decay+}, we get 
 $$ \int_0^T  \|u_t\|_{L^{3,1}}\,dt\leq K_0,$$
 which  completes the proof of \eqref{eq:Lip31}, provided $E_1(0)$ is small enough. 
 \smallbreak
 At this stage, putting together all the estimates we have obtained so far, we deduce
 the decay estimates \eqref{eq:decay}, \eqref{eq:decay+}, and 
  the following control of the solution:
 \begin{multline}\label{eq:final}
 \underset{t\in[0,T]}\sup \Bigl(E_{0}(t)+E_1(t)+E_2(t)+
 \|w(t)\|_{C^{0,1}}\Bigr)
 +\|(u,w)\|_{L^2(0,T;L^\infty)}\\+\|(\nabla u,\nabla w)\|_{L^1(0,T;L^\infty)}
  +\int_0^T\bigl(D_0+\wt D_1+D_2\bigr)dt
 \leq K_0,\end{multline}
  where $K_0$  is a nondecreasing function of  $E_1(0),$ $M_0,$ $\|\rho_0\|_{L^\infty},$ $\|w_0\|_{C^{0,1}}$
   and  $\|u_0\|_{L^{3/2}}.$

 \subsubsection*{Step 7:  More time decay estimates for the energy functionals}
 Here we  show that if, in addition,  $u_0$ belongs to $\dot B^{-3/2}_{2,\infty}$ then   \eqref{eq:energydecaybis} holds.
 The key is to establish that regularity  $\dot B^{-3/2}_{2,\infty}$  of $u$  is conserved through the evolution.
   This is based on some modification 
   of the parabolic maximal regularity estimates in Besov spaces presented before, that has been    
   discovered by Chemin in \cite{Chemin}
  (see also    \cite[Ineq. (3.39)]{BCD}), namely
 \begin{multline}\label{eq:v32}\|u\|_{L^\infty(0,t;\dot B^{-3/2}_{2,\infty})} + \|u\|_{\wt L^1(0,t;\dot B^{1/2}_{2,\infty})}\\\leq \|u_0\|_{\dot B^{-3/2}_{2,\infty}} +  \|u\cdot\nabla u\|_{\wt L^1(0,t;\dot B^{-3/2}_{2,\infty})}
+\|\rho(w-u)\|_{\wt L^1(0,t;\dot B^{-3/2}_{2,\infty})}.\end{multline} 
Above, the notation $\|\cdot\|_{\wt L^1(0,t;\dot B^{1/2}_{2,\infty})}$ stands for some
norm for functions on $(0,t)\times\R^3$ (see \cite[Def. 2.67]{BCD}).
 There is no real need here to recall the definition: we just  have to know that 
 spaces $\wt L^r(0,t;\dot B^\sigma_{2,\infty})$ are closely related to $L^r(0,t;\dot B^\sigma_{2,\infty})$
 and satisfy
$$\|z\|_{\wt L^1(0,t;\dot B^{\sigma}_{2,\infty})}\leq \|z\|_{L^1(0,t;\dot B^{\sigma}_{2,\infty})}.$$
Now, owing to \eqref{eq:embed0}  and  using real interpolation, we have
$$\begin{aligned}\|u\cdot\nabla u\|_{L^1(0,t;\dot B^{-3/2}_{2,\infty})}&\lesssim \|u\cdot\nabla u\|_{L^1(0,t;L^1)}\\
&\leq \|u\|_{L^2(0,t;L^2)}\|\nabla u\|_{L^2(0,t;L^2)}\\&
\lesssim\|u\|_{L^\infty(0,t;\dot B^{-3/2}_{2,\infty})}^{1/3}\|u\|_{\wt L^1(0,t;\dot B^{1/2}_{2,\infty})}^{1/3}
\|\nabla u\|_{L^2(0,t;L^2)}^{4/3}\\&\leq \|u\|_{L^\infty(0,t;\dot B^{-3/2}_{2,\infty})\cap \wt L^1(0,t;\dot B^{1/2}_{2,\infty})}^{2/3}
E_0^{2/3}(0).\end{aligned}$$
We also have, using \eqref{eq:embed0} then \eqref{eq:decay+} for some $\beta\in(1,3/2)$ and Cauchy-Schwarz inequality 
$$\begin{aligned}\|\rho(w-u)\|_{L^1(0,t;\dot B^{-3/2}_{2,\infty})}&\lesssim \|\rho(w-u)\|_{L^1(0,t;L^1)}\\
&\lesssim\|\sqrt\rho\|_{L^\infty(0,t;L^2)} \int_0^t (1+a_1\tau)^{\beta/2}\|\sqrt\rho\,(w-u)\|_{L^2}\: (1+a_1\tau)^{ -\beta/2}\,d\tau\\
&\lesssim e^{CE_0(0)}\sqrt{M_0\,E_1(0)}.\end{aligned}$$
Hence, reverting to \eqref{eq:v32} gives for some absolute constant $C$ and all $t>0,$
\begin{equation}\label{eq:v320}
\|u\|_{L^\infty(0,t;\dot B^{-3/2}_{2,\infty})} + \|u\|_{\wt L^1(0,t;\dot B^{1/2}_{2,\infty})}\leq
C_{0,0}:=C\Bigl(\|u_0\|_{\dot B^{-3/2}_{2,\infty}} + E_0^2(0)+e^{CE_0(0)}\sqrt{M_0E_1(0)}\Bigr)\cdotp
\end{equation}
Starting from \eqref{eq:energybalance}, it is  now possible to improve the time decay of $E_0$ by arguing essentially as we did before for $E_1$: 
using \eqref{eq:interpo} with $z=u,$ then \eqref{eq:v320}, we have
\begin{equation}\label{eq:lb1}
\|\nabla u\|_{L^2}^2\geq C_{0,0}^{-4/3}\bigl(\|u\|_{L^2}^2\bigr)^{ 5/3}.\end{equation}
In order to bound from below the second term of $E_0,$ we use 
the triangle inequality then Sobolev inequality as follows:
$$\begin{aligned}
\|\sqrt \rho\,w\|_{L^2}^2&\leq 2\|\sqrt\rho\,(w-u)\|_{L^2}^2 +2\|\sqrt\rho\, u\|_{L^2}^2\\
&\leq 2\|\sqrt\rho\,(w-u)\|_{L^2}^2 +2\|\rho\|_{L^{3/2}}\|u\|_{L^6}^2\\
&\leq 2\|\sqrt\rho\,(w-u)\|_{L^2}^2 +CM_0^{2/3} R^{1/3}\|\nabla u\|_{L^2}^2.\end{aligned}$$
 We thus have for all $\eps>0,$ 
 $$
 D_0\geq \|\nabla u\|_{L^2}^2+\eps\|\sqrt\rho\, w\|_{L^2}^2 - C\eps M_0^{2/3} R^{1/3}\|\nabla u\|_{L^2}^2.$$
 Choosing $\eps=\eta\min(M_0^{-2/3}R^{-1/3},1)$ for some suitably small (absolute) $\eta,$ and then using \eqref{eq:lb1}
 and reverting to \eqref{eq:energybalance}, we discover that
 $$ \frac d{dt}E_0 +C_{0,0}^{-4/3}\bigl(\|u\|_{L^2}^2\bigr)^{ 5/3}+\eps\|\sqrt\rho\, w\|_{L^2}^2\leq	0.$$
 Since $\|\sqrt\rho\, w\|_{L^2}^2\leq E_0(0)$ and $E_0(0)$ has been assumed to be small, 
 the above inequality is still valid  if $\|\sqrt\rho\, w\|_{L^2}^2$
 is raised to the power $5/3,$ giving eventually for some constant $a_0$ depending only on norms of the data, 
 $$\frac d{dt}E_0+a_0E_0^{5 /3}\leq0,$$
 whence the first part of Inequality \eqref{eq:energydecaybis}. 
 Of course, as for $E_1$ in \eqref{eq:decay}, this implies that 
 $$ \int_0^T ( 1+a_0t)^{\beta} D_0(t)\,dt\leq\frac{\beta}{3/2-\beta}\: E_0(0)\quad\hbox{for all }\ \beta\in(0,3/2).$$
  Observe that  Inequality \eqref{eq:v320} also allows to improve the decay of $E_1.$ Indeed, it suffices to combine it
with  the interpolation inequality 
  $$
  \|\nabla u\|_{L^2}\lesssim \|\nabla^2u\|_{L^2}^{5/7} \|u\|_{\dot B^{-3/2}_{2,\infty}}^{2/7}$$
  to get 
  $$  \|\nabla^2u\|_{L^2}^2\geq C_{0,0}^{-4/5}\bigl(\|\nabla u\|_{L^2}^2\bigr)^{7/5}.$$
  Then, we can faithfully follow the argument leading to \eqref{eq:energydecay}, and get the second
  inequality of \eqref{eq:energydecaybis}. 
  
   \subsubsection*{Step 8: Large time behavior of the density}
   We here assume that the solution of (ENS) under consideration is defined for all positive time. 
   To study the convergence of $\rho(t)$ for $t\to\infty,$ the starting point is that
  \begin{equation}\label{eq:rhot}  \rho(t)=\rho_0-\int_0^t\div(\rho w)\,d\tau=\rho_0-\div\int_0^t \rho u\,d\tau-\div\int_0^t\rho(w-u)\,d\tau.\end{equation}
  Now, we have
  $$  \|\rho(w-u)\|_{L^1} \leq \|\sqrt\rho\|_{L^2}\|\sqrt\rho\,(w-u)\|_{L^2}\leq\sqrt{M_0}\sqrt{E_1}\lesssim (1+t)^{-5/4},$$
  and thus $\|\rho(w-u)\|_{L^1}$ is integrable on $\R_+.$ 
Note that   $\|\rho(w-u)\|_{L^2}$ is also integrable on $\R_+,$ since  
   $$  \|\rho(w-u)\|_{L^2} \leq \|\sqrt\rho\|_{L^\infty}\|\sqrt\rho\,(w-u)\|_{L^2}\lesssim (1+t)^{-5/4}.$$
 Similarly, \begin{align}\label{eq:rho1}
    \|\rho u\|_{L^1} &\leq \|\rho\|_{L^{6 /5}}\|u\|_{L^6}\lesssim  M_0^{5/6}\|\rho\|_{L^\infty}^{1/6} \|\nabla u\|_{L^2}
   \lesssim (1+t)^{-5/4}\\\andf \label{eq:rho2}
    \|\rho u\|_{L^2} &\leq \|\rho\|_{L^3}\|u\|_{L^6}\lesssim M_0^{1/3}\|\rho\|_{L^\infty}^{2/3} \|\nabla u\|_{L^2}
   \lesssim  (1+t)^{-5/4}.\end{align}
   Hence $ \|\rho u\|_{L^1}$ and  $\|\rho u\|_{L^2}$ are integrable on $\R_+$  and 
      $$\rho_\infty:=\rho_0-\int_0^\infty\div(\rho w)\,d\tau$$
      is thus well-defined as an element of $\rho_0+\bigl(\dot W^{-1,1}\cap \dot H^{-1}\bigr).$ 
Remembering \eqref{eq:rhot}, we get 
   $$\rho_\infty-\rho(t)=-\int_t^\infty\div(\rho w),$$
   whence    \eqref{eq:convdensity}, thanks to \eqref{eq:rho1} and \eqref{eq:rho2}.
   Finally, the fact that $\rho(t)$ is uniformly bounded in $L^\infty$ ensures that $\rho_\infty$ is bounded, too. \qed


 \section{The proof of Theorems \ref{thm:main1} and \ref{thm:main1b}} \label{s:proofs}

 In this section, we take advantage of the a priori estimates of the previous section to prove our global existence and 
 uniqueness statements pertaining to the case of initial velocity in $H^1.$

 \subsection{Uniqueness and stability} \label{s:uniqueness}

 We consider two solutions $(\rho_i,w_i,u_i,P_i)$ given by Theorem \ref{thm:main1},
 and set 
 $$
 \dr:=\rho_2-\rho_1,\quad \dw:=w_2-w_1,\quad \du:=u_2-u_1\andf \dP:=P_2-P_1.
 $$
The equations satisfied  by $(\rho_2,\dw,\du,\dP)$ may be formulated
as a linear system of type \eqref{eq:LENS}, namely, 
\begin{equation}\label{eq:uniq}
 \left\{\begin{aligned}
&\rho_{2,t}+\div(\rho_2 w_2)=0,\\
&\dw_t+w_2\cdot\nabla\dw+\dw-\du=-\dw\cdot\nabla w_1,\\
&\du_t+u_2\cdot\nabla\du-\Delta \du+\nabla \dP+\rho_2(\du-\dw)= \dr(w_1-u_1)-\du\cdot\nabla u_1,\\
&\div \du=0.
\end{aligned}\right.
\end{equation}
Hence, taking advantage of \eqref{eq:energybasic} gives us: 
 \begin{multline}\label{eq:uniq1}
 \frac12\frac d{dt}\dE+\dD=-\int \rho_2(\dw\cdot\nabla w_1)\cdot\dw\,dx
 +\int \dr (w_1-u_1)\cdot\du\,dx-\int  (\du\cdot\nabla u_1)\cdot \du\,dx\\\with
\dE:=  \|\sqrt{\rho_2}\dw\|_{L^2}^2+\|\du\|_{L^2}^2\andf
\dD:= \|\nabla\du\|_{L^2}^2 +\|\sqrt{\rho_2}\,(\dw-\du)\|_{L^2}^2.
\end{multline}
Using H\"older, Sobolev and Young inequalities allows to bound the right-hand side as follows:
$$
\begin{aligned}
-\int \rho_2(\dw\cdot\nabla w_1)\cdot\dw\,dx&\leq \|\nabla w_1\|_{L^\infty}\|\sqrt{\rho_2}\,\dw\|_{L^2}^2,\\
 -\int  (\du\cdot\nabla u_1)\cdot \du\,dx&\leq \|\nabla u_1\|_{L^\infty}\|\du\|_{L^2}^2,\\
 \int \dr (w_1-u_1)\cdot\du\,dx&\leq \|\dr\|_{\dot H^{-1}} \|\nabla( (w_1-u_1)\cdot\du)\|_{L^2}\\
 &\leq \|\dr\|_{\dot H^{-1}} \bigl( \|\nabla(w_1-u_1)\|_{L^\infty}\|\du\|_{L^2}+\|w_1-u_1\|_{L^\infty}\|\nabla\du\|_{L^2}\bigr)\\
 &\leq  \frac12\Bigl(\|\dr\|_{\dot H^{-1}}+\|\du\|_{L^2}^2\Bigr) \|\nabla(w_1-u_1)\|_{L^\infty}\\&\hspace{4cm}+
 \frac12\|\nabla\du\|_{L^2}^2+\frac12\|w_1-u_1\|_{L^\infty}^2\|\dr\|_{\dot H^{-1}}^2.
\end{aligned}$$
Hence, reverting to \eqref{eq:uniq1}, we get 
 \begin{equation}\label{eq:uniq2}
 \frac d{dt}\dE+\dD\leq 3\|(\nabla w_1,\nabla u_1)\|_{L^\infty}\dE
 + \bigl(\|w_1-u_1\|_{L^\infty}^2+\|\nabla(w_1-u_1)\|_{L^\infty}\bigr)\|\dr\|_{\dot H^{-1}}^2.
 \end{equation}
To bound $\|\dr\|_{\dot H^{-1}}^2,$  we argue by duality. Let us fix some $T>0$ and write that
$$\|\dr(T)\|_{\dot H^{-1}}=\underset{\|\phi\|_{\dot H^1}=1}\sup\int \phi\,\dr(T)\,dx.$$
 Let us take some $\phi$ in the unit sphere of $\dot H^1,$ and  solve the following backward transport equation:
 \begin{equation}\label{eq:backward} \d_t r+w_1\cdot\nabla r=0,\qquad r|_{t=T}=\phi.\end{equation}
 From  \eqref{eq:backward} and integration by parts, we get 
 $$\begin{aligned}
 0&=\int_0^T\!\!\int_{\R^3}\dr\bigl(\d_tr+w_1\cdot\nabla r)\,dx\,dt\\
 &=-\int_0^T\!\!\int_{\R^3}\bigl(\d_t\dr+\div(w_1\dr)\bigr)r\,dx\,dt
 +\int_{\R^3} \phi\,\dr(T)\,dx-\int_{\R^3} r(0)\,\dr(0)\,dx.\end{aligned}$$
Hence, using the fact that 
 $$ \d_t\dr+\div(w_1\dr)=-\div(\rho_2\dw),$$ we discover that 
$$\begin{aligned}\int_{\R^3} \phi\,\dr(T)\,dx&=\int_{\R^3} r(0)\,\dr(0)\,dx-\int_0^T\!\!\int_{\R^3} r\div(\rho_2\dw)\,dx\,dt\\
&\leq \|\dr(0)\|_{\dot H^{-1}}\|\nabla r(0)\|_{L^2}+\int_0^T\|\rho_2\dw\|_{L^2}\|\nabla r\|_{L^2}\,dt.\end{aligned}$$
Now, from  $\|\nabla\phi\|_{L^2}=1$ and classical estimates for the transport equation, we infer that
$$\|\nabla r(t)\|_{L^2}\leq e^{\int_t^T\|\nabla w_1\|_{L^\infty}\,d\tau},\quad t\in[0,T].$$
Hence, we end up with the following bound for $\dr$:
\begin{equation}\label{eq:uniq4}
\|\dr(T)\|_{\dot H^{-1}}\leq e^{\int_0^T\|\nabla w_1\|_{L^\infty}\,dt}\biggl(\|\dr(0)\|_{\dot H^{-1}}+
\int_0^Te^{-\int_0^t\|\nabla w_1\|_{L^\infty}\,d\tau}\|\rho_2\dw\|_{L^2}\,dt\biggr)\cdotp\end{equation}
Applying Gronwall Lemma to \eqref{eq:uniq2},  
then  the above inequality yields
\begin{multline*}
\dE(t)+\int_0^t\dD\,d\tau\leq  e^{3\int_0^t\|(\nabla w_1,\nabla u_1)\|_{L^\infty}\,d\tau}
\biggl(\dE(0)\\
+\int_0^t  \bigl(\|w_1-u_1\|_{L^\infty}^2+\|\nabla(w_1-u_1)\|_{L^\infty}\bigr)
\biggl(\|\dr(0)\|_{\dot H^{-1}}^2 +\biggl(\int_0^\tau
\|\rho_2\dw\|_{L^2}\,d\tau'\biggr)^2\biggr)d\tau\biggr)\cdotp
\end{multline*}
By Cauchy-Schwarz inequality, we have
$$\biggl(\int_0^\tau
\|\rho_2\dw\|_{L^2}\,d\tau\biggr)^2\leq \tau \int_0^\tau\|\rho_2\dw\|_{L^2}^2\,d\tau.$$
Finally, using again Gronwall lemma, and the fact that $(w_1-u_1)\in L^2(\R_+;L^\infty)$ 
and that $(\nabla w_1,\nabla u_1)\in L^1(\R_+;L^\infty),$ we get
\begin{multline}\label{eq:stability}\underset{\tau\in[0,t]}\sup\dE(\tau)+\int_0^t\dD\,d\tau\leq K(t)\biggl(\dE(0)
\\+\biggl(\int_0^t\tau \bigl(\|w_1-u_1\|_{L^\infty}^2+\|\nabla(w_1-u_1)\|_{L^\infty}\bigr)\,d\tau\biggr)
\|\dr(0)\|_{\dot H^{-1}}^2\biggr)\end{multline}
for some bounded positive function $K,$
from which one can conclude to uniqueness, in the case where the solutions 
$(\rho_1,w_1,u_1,P_1)$ and $(\rho_2,w_2,u_2,P_2)$ originate from  the same  data.\qed 


 \subsection{The existence scheme}
 
 The overall strategy is  standard : we smooth out the data, and obtain a sequence of local-in-time regular solutions. 
 Then, leveraging a blow-up criterion involving the Lipschitz norms of the velocities, and combining with the Lipschitz control 
 that has been obtained in the previous section, we show that these regular solutions are actually global. 
In particular, they satisfy \eqref{eq:final}, which provides us with uniform estimates in the solution space described in 
Theorem \ref{thm:main1}. Furthermore, all the decay estimates that we proved so far  are satisfied uniformly by 
the sequence of smooth solutions.
The final step is to prove convergence : it is based both on the stability estimates of the previous subsection
(which ensures the strong convergence of the Navier-Stokes velocity and of the density), and on a compactness argument (to handle the Euler velocity).

\subsubsection*{Step 1: construction of a sequence of local-in-time regular solution for approximated data}

 Take a nonnegative function $\chi$ in $\cC_0^\infty(B(0,1))$ such that $\|\chi\|_{L^1}=1$ and set
$\chi^n(x):=n^{-3}\chi(n^{-1}x)$ for all $n\in\N^*.$ 
Fix some initial data $(\rho_0,w_0,u_0)$ satisfying the assumptions of Theorem \ref{thm:main1} and consider System (ENS) with initial data
$$(\rho_0^n,w_0^n,u_0^n):=(\chi^n\star \rho_0, \chi^n\star w_0,\chi^n \star u_0).$$
Clearly, the triplet $(\rho_0^n,w_0^n,u_0^n)$ belongs to all Sobolev spaces.
Now, since (ENS) is a coupling between a transport equation, a (damped) Burgers equation with source term, 
and a parabolic equation with a nonlocal quadratic nonlinearity, it is easy to establish 
that is has a unique local-in-time solution $(\rho^n,w^n,u^n,\nabla P^n)$ belonging  to all Sobolev spaces (one can use for instance the Friedrichs mollification method
described in \cite[Chap. 4 and 9]{BCD}). 

For the sequel, it suffices  to consider just the Sobolev regularity $H^2$ for  $\rho^n,$ and $H^3$ for $(w^n,u^n).$ We
omit talking about the pressure, since it may be computed from the third equation of $(ENS).$ 
We denote by $T^n$ the supremum of all positive  times $T$ such that $(\rho^n,w^n,u^n)$ belongs to 
$$ E_T:=\Bigl\{(\rho,w,u)\in \cC([0,T];H^2))\times \cC([0,T];H^3)\times \bigl(\cC([0,T];H^3)\cap  L^2(0,T;H^4)\bigr)\Bigr\}\cdotp$$ 
 Note that since $\rho_0\in L^1$ and $u_0\in L^{3/2},$ we also have 
 $$ \rho^n\in\cC([0,T^n);L^1)\andf u^n\in\cC([0,T^n);L^{3/2}).$$

\subsubsection*{Step 2: A continuation criterion}

By using a basic energy method and classical estimates in Sobolev spaces, it is easy to get  for all $t\in[0,T^n),$ 
\begin{align}\label{est:1}
&\|\rho^n(t)\|_{H^2} \leq \|\rho^n_0\|_{H^2}+C\int_0^t\|\nabla w^n\|_{L^\infty} \|\rho^n\|_{H^2}\,d\tau
+C\int_0^t\|\rho^n\|_{L^\infty} \|w^n\|_{H^3}\,d\tau,\\\label{est:2}
&\|w^n(t)\|_{H^3} + \int_0^t\|w^n\|_{H^3}\,d\tau\leq \|w^n_0\|_{H^3} +  \int_0^t\|u^n\|_{H^3}\,d\tau
+C\int_0^t\|\nabla w^n\|_{L^\infty} \|w^n\|_{H^3}\,d\tau,\\ 
&\|u^n(t)\|_{H^3} + \|u^n\|_{L^2(0,t;H^4)}\leq \|u^n_0\|_{H^3} 
+C\int_0^t\|\nabla u^n\|_{L^\infty} \|u^n\|_{H^3}\,d\tau\label{est:3}\\\nonumber 
&\hspace{3cm} +C \|\rho^n\|_{L^\infty(0,t;L^\infty)}\|w^n\|_{L^2(0,t;H^2)}
 + C\|\rho^n\|_{L^\infty(0,t;H^2)}\|w^n\|_{L^2(0,t;L^\infty)}.\end{align}
Using the Sobolev embedding $H^s\hookrightarrow L^\infty$ for $s>3/2$ and denoting 
$$E^n(t):= \underset{\tau\in[0,t]}\sup \|(\rho^n,w^n,u^n)\|_{H^3} + \|u^n\|_{L^2(0,t;H^4)},$$
we deduce  from the above inequalities that for all $t\in[0,T^n),$
$$E^n(t)\leq  E^n(0) + Ct\bigl(E^n(t)+(E^n(t))^2\bigr) +\sqrt t (E^n(t))^2.$$
From this inequality and a bootstrap argument, we can conclude that there exists some universal constant $c$ such that
$$T^n\geq c\min\Bigl(1,\bigl(\|\rho^n_0\|_{H^2}+\|w^n_0\|_{H^3}+\|u^n_0\|_{H^3}\bigr)^{-2}\Bigr)\cdotp$$
By classical argument, this entails that if $T^n$ is finite, then
$$\underset{t\to T^n}\lim \bigl(\|\rho^n(t)\|_{H^2}+\|w^n(t)\|_{H^3}+\|u^n(t)\|_{H^3}\bigr)=\infty.$$
We claim that  whenever
\begin{equation}\label{eq:blowup2}
\|(\nabla w^n,\nabla u^n)\|_{L^1(0,T;L^\infty)}+\|w^n\|_{L^2(0,T;L^\infty)}+\|\rho^n\|_{L^\infty(0,T;L^\infty)}<\infty,
\end{equation}
the solution may be continued beyond $T.$
\medbreak
Indeed,  assume that \eqref{eq:blowup2} is satisfied for some $T=T^n,$ and
set 
$$R:=\|w^n\|_{L^2(0,T;L^\infty)}+\|\rho^n\|_{L^\infty(0,T;L^\infty)}.
$$
Let us take $c>0$ small enough and set 
$$E^n_R(t):= \|\rho^n(t)\|_{H^2}+\|w^n(t)\|_{H^3}+ \|w^n\|_{L^1(0,t;H^3)}+cR^{-1}\|u^n(t)\|_{H^3} 
+cR^{-1}\|u^n\|_{L^2(0,t;H^4)}. $$
Then,  combining \eqref{est:1}, \eqref{est:2} and \eqref{est:3} yields  for all $t\in[0,T],$
\begin{equation*} E^n_R(t)
\leq  E^n_R(0) + C\int_0^t\bigl(1+\|(\nabla w^n,\nabla u^n)\|_{L^\infty}\bigr) E^n_R\,d\tau
 	\end{equation*}
	and using Gronwall lemma eventually gives
	$$	\underset{t\in[0,T]}\sup \bigl(\|\rho^n(t)\|_{H^2}+\|w^n(t)\|_{H^3}+\|u^n(t)\|_{H^3}\bigr)<\infty,$$
whence our claim.

\subsubsection*{Step 3: Global existence of a smooth solution for regularized data}

Since the constructed solution is rather smooth, it satisfies 
the basic conservation laws \eqref{eq:M0} and \eqref{eq:energybalance} as well as
all the a priori estimates of the previous section on $[0,T^n),$ in particular
\eqref{eq:final} provided $(\rho^n_0,w^n_0,u^n_0)$ satisfies the smallness condition \eqref{eq:smalldata}. 
Since this latter condition holds if it is satisfied by $(\rho_0,w_0,u_0),$ we conclude
that we have  for all $T\in[0,T^n)$ (with obvious notation):
\begin{multline}\label{eq:finaln}
 \underset{t\in[0,T]}\sup \Bigl(E_{0}^n(t)+E_1^n(t)+E_2^n(t)+
 \|w^n(t)\|_{C^{0,1}}\Bigr)
 +\|(u^n,w^n)\|_{L^2(0,T;L^\infty)}\\+\|(\nabla u^n,\nabla w^n)\|_{L^1(0,T;L^\infty)}
  +\int_0^T\bigl(D_0^n+\wt D_1^n+D_2^n\bigr)dt
 \leq K_0,\end{multline}
  where $K_0$  depends only on  $E_1(0),$ $M_0,$ $\|\rho_0\|_{L^\infty},$ $\|w_0\|_{C^{0,1}}$
   and  $\|u_0\|_{L^{3/2}}.$ 
\medbreak
In particular, we have \eqref{eq:blowup2} with $T=T^n,$ and thus $T^n$ must be infinite.
In other words, the solution global. 
Furthermore,  all the norms of  $(\rho^n,w^n,u^n)$ coming into play 
in \eqref{eq:finaln} are  uniformly bounded on $\R_+,$ independently of $n.$

\subsubsection*{Step 4: Passing to the limit}

The results of the previous step readily ensure that, up to subsequence, 
\begin{itemize}
\item $\rho^n\rightharpoonup \rho$ in $L^\infty(\R_+\times\R^3)$ weak *;
\item $w^n\rightharpoonup w$ in $L^\infty(\R_+\times\R^3)$ weak *;
\item $u^n\rightharpoonup u$ in $L^\infty(\R_+;H^1)$ weak *.
\end{itemize}
Morever, the stability estimate \eqref{eq:stability} guarantees that 
$(u^n)$ (resp. $\rho^n$) is a Cauchy sequence in $L^\infty_{loc}(\R_+;L^2)$ (resp. $L^\infty_{loc}(\R_+;\dot H^{-1})$).
Hence $u^n\to u$ and $\rho^n\to \rho$ strongly in these spaces. 
\smallbreak
To upgrade the convergence of $(w^n)_{n\in\N},$ we can notice that, since
$$
w^n_t=-w^n\cdot\nabla w^n+(u^n-w^n),$$
Inequality \eqref{eq:finaln} ensures that 
$(w^n_t)_{n\in\N}$ is bounded in e.g. $L^2(\R_+;L^\infty).$ Since this inequality also ensures that
$(w^n)_{n\in\N}$ is bounded in $L^\infty(\R_+;C^{0,1}),$ taking advantage of Arzela-Ascoli's theorem 
and of Cantor diagonal process allows to conclude that
$w^n\to w$ uniformly on any compact set of $[0,\infty)\times\R^3.$
This allows to pass to the limit in all the terms of System (ENS), and to conclude that
$(\rho,w,u)$ is indeed a solution.

The remaining points (like e.g. time continuity) are standard, and are thus omitted. \qed


\section{The case with critical regularity}\label{s:critical}

Here we consider the case where $u_0$ is only in the Besov space  $B^{1/2}_{2,1}$ instead of the Sobolev space $H^1.$
They key idea is to establish that if $c_0$ in \eqref{eq:smalldata2} is small enough, then 
 there exists some $t_0\in(0,1)$ such that all the conditions of Theorem \ref{thm:main1} are satisfied at time $t=t_0.$
 
 For conciseness,  we will  just derive the key a priori estimates leading to this result, since the proof of short time existence
 in this regularity context is very similar to that of the previous section.  Throughout this section, we will  be using more properties of Besov spaces than in the previous section. 
 The reader is invited to consult \cite[Chap.2 ]{BCD} for more details. 
 \medbreak
 To proceed, we assume that we are given a solution on the time interval $[0,T]$ satisfying suitable regularity properties and such that
\begin{equation}\label{eq:bootstrap}
U:=\|u\|_{L^2(0,T;L^\infty)}<\infty\andf \|\nabla w\|_{L^1(0,T;L^\infty)}\leq \log2.\end{equation}
 Let us first observe that \eqref{eq:rhoLinfty}  and \eqref{eq:bootstrap}  immediately imply that
\begin{equation}\label{eq:rho}\|\rho\|_{L^\infty(0,T;L^\infty)} \leq 2\|\rho_0\|_{L^\infty}.\end{equation}
Let us define 
$$
u_L^1(t):= e^{t\Delta}u_0,\quad u_L^2(t):=\int_0^te^{(t-\tau)\Delta} (\rho(w-u))(\tau)\,d\tau
\andf \wt u:=u-u_L^1-u_L^2.$$
From parabolic maximal regularity (see \cite{BCD,Chemin}), we get 
\begin{align}\label{eq:vL1}
\|u_L^1\|_{L^\infty(0,T;\dot B^{1/2}_{2,1})}+\|u_L^1\|_{L^1(0,T;\dot B^{5/2}_{2,1})} &\leq C\|u_0\|_{\dot B^{1/2}_{2,1}},\\[1ex]
\|u_L^2\|_{L^\infty(0,T;\dot H^{1/2})}+\|u_L^2\|_{L^{4/3}(0,T;\dot H^2)} &\leq C\|\rho(w-u)\|_{L^{4/3}(0,T;L^2)}\nonumber\\
&\leq CT^{1/4} \|\rho\|_{L^\infty(0,T;L^\infty)}^{1/2} \|\sqrt\rho\,(w-u)\|_{L^2(0,T;L^2)}.\nonumber\end{align}
Combining this latter inequality with interpolation, \eqref{eq:energybasic} and \eqref{eq:rho} gives 
\begin{equation}\label{eq:vL2}
\|u_L^2\|_{L^2(0,T;\dot B^{3/2}_{2,1})}\lesssim  T^{1/4}\sqrt{ \|\rho_0\|_{L^\infty} E_0(0)}\cdotp
\end{equation}
The fluctuation $\wt u$ satisfies the evolutionary Stokes equation
$$\wt u_t-\Delta\wt u+\nabla\wt P=-u\cdot\nabla u,\qquad \wt u|_{t=0}=0.$$
Let   $u_L:=u_L^1+u_L^2.$ Using once more parabolic maximal regularity in Besov spaces,
then, the facts that the projector $\cP$ on divergence-free vector fields  is continuous on $ \dot B^{1/2}_{2,1}$
and that the pointwise product maps $\dot B^{1/2}_{2,1}\times \dot B^{3/2}_{2,1}$ to   $\dot B^{1/2}_{2,1},$ we get 
$$\begin{aligned}
\|\wt u\|_{L^\infty(0,T;\dot B^{1/2}_{2,1})}+\|\wt u\|_{L^2(0,T;\dot B^{3/2}_{2,1})}+
\|\wt u\|_{L^1(0,T;\dot B^{5/2}_{2,1})}&\lesssim \|\cP(u\cdot\nabla u)\|_{L^1(0,T;\dot B^{1/2}_{2,1})}\\
&\lesssim \|u\cdot\nabla u\|_{L^1(0,T;\dot B^{1/2}_{2,1})}\\
&\lesssim  \|\wt u\|_{L^2(0,T;\dot B^{3/2}_{2,1})}^2+  \|u_L\|_{L^2(0,T;\dot B^{3/2}_{2,1})}^2.
\end{aligned}$$
Consequently, provided 
  $\|u_L\|_{L^2(0,T;\dot B^{3/2}_{2,1})}$ is small enough, we have 
  the following control on $\wt u,$ due to \eqref{eq:vL1} and \eqref{eq:vL2}: 
  \begin{equation}\label{eq:wtv}
  \|\wt u\|_{L^\infty(0,T;\dot B^{1/2}_{2,1})}+\|\wt u\|_{L^1(0,T;\dot B^{5/2}_{2,1})}\lesssim 
 \|u_0\|_{\dot B^{1/2}_{2,1}}^2+ T^{1/2} \|\rho_0\|_{L^\infty} E_0(0). 
\end{equation}
Owing to the embedding 
\begin{equation}\label{eq:embedd} \dot B^{3/2}_{2,1}\hookrightarrow L^\infty \end{equation}  and to Inequalities \eqref{eq:vL1}, 
\eqref{eq:vL2} and \eqref{eq:wtv}, we thus have 
\begin{equation}\label{eq:vL2Linfty}
U\leq C\bigl(  \|u_0\|_{\dot B^{1/2}_{2,1}}+\|u_0\|_{\dot B^{1/2}_{2,1}}^2+
 T^{1/4} \|\rho_0\|_{L^\infty}^{1/2} E_0(0)^{1/2}+ T^{1/2} \|\rho_0\|_{L^\infty} E_0(0)\bigr)\cdotp
 \end{equation}
Note that   \eqref{eq:vL1}, \eqref{eq:wtv} and \eqref{eq:embedd}  provide us with  a  Lipschitz control  of  $u_L^1$ and $\wt u.$ 
To achieve a Lipschitz control on $u_L^2,$ one can  use \eqref{eq:embed} and the following  maximal regularity estimate
in Lorentz spaces for the heat equation\footnote{Owing to the definition of  $L^{3,1}$ in \eqref{def:L31}, this inequality 
 may be deduced from the corresponding ones in Lebesgue spaces $L^2(\R^3)$ and $L^6(\R^3)$, and real 
interpolation.}:
 $$\|\nabla^2 u_L^2\|_{L^2(0,T;L^{3,1})} \lesssim \|\rho(w-u)\|_{L^2(0,T;L^{3,1})}.$$
Hence, remembering \eqref{eq:rho} and using interpolation and  H\"older inequalities, 
\begin{equation}\label{eq:nablav}
\|\nabla u_L^2\|_{L^2(0,T;L^\infty)} \leq \|\rho\|_{L^\infty(0,T;L^\infty)}^{2/3}\|\sqrt\rho\,(w-u)\|_{L^2(0,T;L^2)}^{2/3}
\|(w,u)\|_{L^2(0,T;L^\infty)}^{1/3}.\end{equation}
Taking advantage of  \eqref{eq:uL2} and of  \eqref{eq:bootstrap}, we get
\begin{equation}\label{eq:uu}\|w\|_{L^2(0,T;L^\infty)}\leq 2\|w_0\|_{L^\infty}+ 2U.\end{equation}
Hence, reverting to \eqref{eq:nablav}, we find that
$$\|\nabla u_L^2\|_{L^2(0,T;L^\infty)} \lesssim  \|\rho_0\|_{L^\infty}^{2/3} E_0^{1/3}(0)\bigl(\|w_0\|_{L^\infty}^{1/3} +U^{1/3}\bigr),$$
and thus,  if $\|u_0\|_{\dot B^{1/2}_{2,1}}$ is small enough, 
$$\|\nabla u\|_{L^1(0,T;L^\infty)} \lesssim 
 \|u_0\|_{\dot B^{1/2}_{2,1}}+ T^{1/2} \bigl(\|\rho_0\|_{L^\infty} E_0(0) +  \|\rho_0\|_{L^\infty}^{2/3} E_0^{1/3}(0)\bigl(\|w_0\|_{L^\infty}^{1/3} +U^{1/3}\bigr)\bigr)\cdotp$$
 Using also \eqref{eq:duinfty} gives
 $$\|\nabla w\|_{L^1(0,T;L^\infty)}\leq 2\|\nabla w_0\|_{L^\infty}+ 2\|\nabla u\|_{L^1(0,T;L^\infty)},$$
 whence
 \begin{multline*} 
\|\nabla w\|_{L^1(0,T;L^\infty)} \lesssim \|\nabla w_0\|_{L^\infty}+
 \|u_0\|_{\dot B^{1/2}_{2,1}}+ T^{1/2} \bigl(\|\rho_0\|_{L^\infty} E_0(0)\\ +  \|\rho_0\|_{L^\infty}^{2/3} E_0^{1/3}(0)\bigl(\|w_0\|_{L^\infty}^{1/3} +U^{1/3}\bigr)\bigr)\cdotp
\end{multline*}
Putting together with \eqref{eq:vL2Linfty}, one can conclude that one can find some small enough constant $c'_0$
such that if 
\begin{multline}\label{eq:c'0} \|\nabla w_0\|_{L^\infty}+
 \|u_0\|_{\dot B^{1/2}_{2,1}}+\|\rho_0\|_{L^\infty} E_0(0) \\+  \|\rho_0\|_{L^\infty}^{2/3} E_0^{1/3}(0)
 \bigl(\|w_0\|_{L^\infty}^{1/3}+\|u_0\|_{\dot B^{1/2}_{2,1}}^{1/3}+\|\rho\|_{L^\infty}^{1/6} E_0^{1/6}(0)\bigr)\leq c'_0,\end{multline} then, 
 the second inequality of \eqref{eq:bootstrap} does  hold with strict inequality provided
 $T\leq1.$ Combining with a local-in-time existence result (the proof of which is left to the reader)
 and a bootstrap argument, we conclude that if \eqref{eq:c'0} holds then  existence holds
 true on $[0,1].$
 \smallbreak
 In particular, the energy balance \eqref{eq:energybalance} holds  and one can thus find some $t_0\in[0,1]$ such that  $E_{1}(t_0)$ satisfies the smallness condition of 
 Theorem \ref{thm:main1}. 
 Indeed,  that $\|\nabla w(t_0)\|_{L^\infty}$ is small for all $t_0\in[0,1]$ stems from \eqref{eq:smalldata2} and \eqref{eq:dwLinfty}
 and the basic energy balance \eqref{eq:energybalance} with $t=1$ implies that there exists $t_0\in[0,1]$ such that
 $$ \|\nabla u(t_0)\|_{L^2}^2 +\|\sqrt{\rho(t_0)}\,(w-u)(t_0)\|_{L^2}^2 \leq E_0(0).$$
 Hence, if $E_0(0)$ is small enough, then  the smallness condition \eqref{eq:smalldata2} is satisfied  at $t=t_0.$ 
 From that point, one can use the uniqueness result to continue the solution into a $H^1$ type solution, beyond time $t=t_0,$
 that satisfies the properties of Theorem \ref{thm:main2}.

   \begin{small}	 

\end{small}


\begin{thebibliography}{99}

\bibitem{AB} O. Anoshchenko and A. Boutet de Monvel~:  The existence of the global generalized solution of the system of equations describing suspension motion,
\emph{Math. Methods Appl. Sci.}, {\bf  20} (1997), no. 6, 495--519. 

 \bibitem{BCD} H.  Bahouri,  J.-Y.  Chemin and   R.  Danchin~:  \emph{Fourier  Analysis  and  Nonlinear  Partial  Differential  Equations,}  Grundlehren  der  Mathematischen  Wissenschaften,  vol. {\bf 343},  Springer-Verlag,  Berlin,  Heidelberg,  2011.
 
 \bibitem{BEHK}  A. Baradat, L. Ertzbischoff and D. Han-Kwan~: Multiphasic formulation of Vlasov
equations and applications,  work in progress.
 
   \bibitem{BDGM} L. Boudin, L.  Desvillettes, C.  Grandmont and A. Moussa~: Global existence of solutions for the coupled Vlasov and Navier-Stokes equations, \emph{Differential Integral Equations}, {\bf 22}(11-12),  (2009),  1247--1271. 
 
\bibitem{Chemin} J.-Y. Chemin~: Th\'eor\`emes d'unicit\'e pour le syst\`eme de Navier-Stokes tridimensionnel, \emph{J. Anal. Math.},  {\bf 77} (1999),  27--50. 

\bibitem{CK} Y.-P. Choi and J. Jung ~: On the Cauchy problem for the pressureless Euler-Navier-Stokes system in the whole
space, \emph{J. Math. Fluid Mech.}, {\bf 23}  (2021), no. 4, Paper No. 99


\bibitem{CJK} 
Y.-P. Choi, J. Jung  and J. Ki~ : A revisit to the pressureless Euler-Navier-Stokes system in the whole space and its optimal temporal decay,  \emph{J. Differential Equations}, {\bf 401}  (2024), 231--281.   

\bibitem{R-VNS} R. Danchin~ : Fujita-Kato solutions and optimal time decay  for the Vlasov-Navier-Stokes system
in the whole space, arXiv:2405.09937. 

\bibitem{Ertz}  L.  Ertzbischoff~:  Global derivation of a Boussinesq-Navier-Stokes type system from fluid-kinetic equations,
\emph{Ann. Fac. Sci. Toulouse}, {\bf 33}(4), 1059--1154, 2024. 


\bibitem{FK} H. Fujita and T. Kato :  On the Navier-Stokes initial value problem I, {\it Archive for Rational Mechanics
 and Analysis}, {\bf 16},  (1964), 269-315.
 
   \bibitem{DHK} D. Han-Kwan~: Large time behavior of small-data solutions to the Vlasov-Navier-Stokes system on the whole space, 
\emph{Probability and Mathematical Physics}, {\bf 3}(1), (2022), 35--67. 

\bibitem{HKM} D. Han-Kwan and D. Michel~: On hydrodynamic limits of the Vlasov-Navier-Stokes system,
\emph{Mem. Amer. Math. Soc.}, {\bf 302} (2024), no. 1516, v+115 pp. 

\bibitem{HMM} D. Han-Kwan, A. Moussa and I. Moyano: Large time behavior of the Vlasov-Navier-Stokes system on the torus, 
\emph{Archiv for Rational Mechanics and Analysis}, {\bf 236}(3), (2020), 1273--1323. 

\bibitem{HTZ} F. Huang, H. Tang and W. Zou ~: Global well-posedness and optimal time decay rates of solutions
to the pressureless Euler-Navier-Stokes system,
 \emph{Commun. Math. Anal. Appl.} {\bf  3} (2024), no. 4, 582--623.

\bibitem{HTWZ} F. Huang, H. Tang, G. Wu  and W. Zou~ : Global well-posedness and large-time
 behavior of classical solutions to the Euler-Navier-Stokes system in $\R^3,$  
 \emph{J. Differential Equations}, {\bf 410} (2024), 76--112.
 
 \bibitem{Lem} V. Lemari\'e~ : The pressureless Euler-Navier-Stokes system, arXiv:2505.17577.

  \bibitem{Nash} J. Nash~: Continuity of solutions of parabolic and elliptic equations.
{\em Amer. J. Math.}, {\bf 80}, (1958),  931--954. 


\bibitem{Wi} F.A. Williams~: \emph{Combustion theory}, Second Edition. The Benjamin Cummings Publishing, 1985.


\bibitem{Zhai} X. Zhai, Y. Chen, Y. Li and Y. Zhao~: {\it Optimal well-posedness for the pressureless Euler-Navier-Stokes system}, \emph{J. Math. Phys.}  {\bf 64} (2023), no.~5, Paper No. 051506, 13 pp.

\end{thebibliography}
\end{document}